\documentclass[12pt,draftcls,onecolumn]{IEEEtran}

\ifCLASSINFOpdf
\else
\fi

\usepackage{ifpdf}
\ifpdf
    \usepackage{graphicx}
    \usepackage{epstopdf}
\else
    \usepackage{graphicx}

\fi

\usepackage{mathptmx} 
\usepackage{times} 
\usepackage{amsmath} 
\usepackage{amssymb}  
\usepackage{tabularx}
\usepackage{algorithm}
\usepackage{algorithmic}
\usepackage{enumerate}
\usepackage[dvipsnames,usenames]{color}
\usepackage[colorlinks,%
linkcolor=BrickRed,%
filecolor=BrickRed,%
citecolor=RoyalPurple,%
]{hyperref}
\usepackage[font=small]{caption}
\usepackage{subcaption} 

\usepackage{color}

\usepackage{tabularx}

\title{\LARGE \bf
Feedback Particle Filter on Matrix Lie Groups}

\author{Chi Zhang, Amirhossein Taghvaei and Prashant G. Mehta
\thanks{C. Zhang, A. Taghvaei and P. G. Mehta are with the Coordinated Science Laboratory and
the Department of Mechanical Science and Engineering at the University of
Illinois at Urbana-Champaign (UIUC)
{\tt\small czhang54@illinois.edu; taghvae2@illinois.edu; mehtapg@illinois.edu}}
\thanks{Financial support from the NSF CMMI grants 1334987 and 1462773 is gratefully acknowledged.}%
\thanks{The conference versions of this paper appear in~\cite{CDC2016, ACC2016}.}
}

\def\qed{\hspace*{\fill}~\IEEEQED\par\endtrivlist\unskip}

\def\Re{\mathbb{R}}

\usepackage{eufrak}

\def\Sec#1{Sec.~\ref{#1}}

\def\notes#1{\marginpar{\tiny #1}\typeout{Notes!
Notes!
Notes!
}}
\renewcommand{\notes}[1]{\typeout{notes!}}

\newcommand{\field}[1]{\mathbb{#1}}

\def\Re{\field{R}}

\def\Sec#1{Sec.~\ref{#1}}

\def\clG{{\cal G}}

\def\clL{{\cal L}}
\def\clN{{\cal N}}

\def\clV{{\cal V}}
\def\clW{{\cal W}}

\def\K{{\sf K}} 
\def\k{{\sf k}} 
\def\u{{\sf u}} 
 
\def\id{{I}} 
\def\innov{\mathrm{I}} 
\def\lap{\Delta_{\rho_t}} 
\def\Dt{\Delta t}
\def\dt{\ud t}
\def\ds{\ud s}

\def\dq{\ud q}
\def\dR{\ud R}
\def\dZ{\ud Z}
\def\dI{\ud \innov}
\def\dB{\ud B}
\def\dW{\ud W}

\def\UZ{\mathcal{Z}_t}

\def\Lstar{\mathcal{L}^*}
\def\Cinf{C^{\infty}}
\def\Cinfc{\Cinf_c}

\def\x{{x}} 
\def\X{{X}} 
\def\ip{\alpha} 
\def\ipp{\beta} 
\def\io{j} 
\def\dev{\chi} 
\def\cov{\overline{\Sigma}}

\def\ke{k^{(\epsilon)}}

\def\wken{\tilde{k}^{(\epsilon,N)}}
\def\Ten{T^{(\epsilon,N)}}

\def\grad{\text{grad}}
\def\div{\text{div}}

\def\det{\text{det}}
\def\exp{\text{exp}}

\def\ito{It\^{o}}

\makeatletter

\newcommand{\Rom}[1]{\expandafter\@slowromancap\romannumeral #1@}
\makeatother

\renewcommand{\baselinestretch}{1.28}

\newcommand{\ud}{\,\mathrm{d}}

\def\Expect{{\sf E}}

\newcommand{\ve}[1]{[\,#1\,]^{\vee}}

\newcommand{\lr}[2]{\langle #1, #2 \rangle}
\newcommand{\sk}[1]{[\,#1\,]_{\times}}
\newcommand{\degg}[1]{{#1}^{\circ}}
\newcommand{\inv}[1]{{#1}^{-1}}
\newcommand{\normG}[1]{|#1|_{\clG}}

\newcommand{\bpar}[1]{\big( #1 \big)}
\newcommand{\Bpar}[1]{\Big( #1 \Big)}

\newcommand{\tn}[1]{{#1}_{t_n}}

\DeclareMathOperator{\Tr}{Tr}

\def\Expect{{\sf E}}

\usepackage{ulem}

\newtheorem{theorem}{Theorem}
\newtheorem{corollary}{Corollary}
\newtheorem{proposition}{Proposition}

\newtheorem{definition}{Definition}

\newtheorem{remark}{Remark}
\newtheorem{assumption}{Assumption}

\def\Sec#1{Sec.~\ref{#1}}

\makeatletter
\newsavebox\myboxA
\newsavebox\myboxB
\newlength\mylenA

\newcommand*\xoverline[2][0.75]{%
    \sbox{\myboxA}{$\m@th#2$}%
    \setbox\myboxB\null
    \ht\myboxB=\ht\myboxA%
    \dp\myboxB=\dp\myboxA%
    \wd\myboxB=#1\wd\myboxA
    \sbox\myboxB{$\m@th\overline{\copy\myboxB}$}
    \setlength\mylenA{\the\wd\myboxA}
    \addtolength\mylenA{-\the\wd\myboxB}%
    \ifdim\wd\myboxB<\wd\myboxA%
       \rlap{\hskip 0.5\mylenA\usebox\myboxB}{\usebox\myboxA}%
    \else
        \hskip -0.5\mylenA\rlap{\usebox\myboxA}{\hskip 0.5\mylenA\usebox\myboxB}%
    \fi}
\makeatother

\def\sigmaB{\xoverline{\sigma}_B}
\def\sigmaW{\xoverline{\sigma}_W}

\begin{document}
\normalem
\maketitle

\vspace{-0.5in}
\begin{abstract}

This paper is concerned with the problem of continuous-time nonlinear filtering for 
stochastic processes on a connected matrix Lie group. 
The main contribution of this paper is to derive the feedback particle filter (FPF) 
algorithm for this problem. In its general form, the FPF is shown to provide a coordinate-free 
description of the filter that automatically satisfies the geometric constraints of the manifold. 
The particle dynamics are encapsulated in a Stratonovich stochastic differential 
equation that retains the feedback structure of the original (Euclidean) FPF. 
The implementation of the filter requires a solution of a Poisson equation on the Lie group, 
and two numerical algorithms are described for this purpose.
As an example, the FPF is applied to the problem of attitude estimation\,--\,a nonlinear filtering problem on 
the Lie group $SO(3)$. The formulae of the filter are described using both the rotation matrix and the quaternion coordinates. 
Comparisons are also provided between the FPF and some popular algorithms for attitude estimation, 
namely the multiplicative EKF, the unscented quaternion estimator, the left invariant EKF, 
and the invariant ensemble Kalman filter. 
Numerical simulations are presented to illustrate the comparisons. 

\end{abstract}

\section{Introduction}
\label{sec:intro}

There is a growing interest in the nonlinear filtering community to develop geometric approaches 
for handling constrained systems. In many cases, the constraints are described by 
smooth Riemannian manifolds, in particular the matrix Lie groups. Engineering applications 
of filtering on matrix Lie groups include: 
i) attitude estimation of aircrafts \cite{hua2014implementation, barrau2015TAC}; 
ii) visual tracking of humans and objects \cite{kwon2014pf, choi2012robust}; and 
iii) localization of mobile robots \cite{barczyk2015invariant, hesch2013camera}. 
In these applications, the matrix Lie groups of interest include  
the special orthogonal group $SO(3)$ and the special Euclidean group $SE(3)$. 

\subsection{Problem Statement}

We consider the following continuous-time system evolving on a matrix Lie group $G$, 
\begin{subequations}
\begin{align}
 \ud X_t & = X_t \, V_0(X_t) \, \dt + X_t \, V_\ip(X_t) \circ\ud B_t^\ip, \label{eq:dyn} \\
 \ud Z_t & = h(X_t) \, \dt + \ud W_{t}, \label{eq:obs}
\end{align}
\end{subequations}
where $X_t\in G$ is the state at time $t$, $Z_t\in\Re^m$ is the observation vector; 
$V_0:G\rightarrow\clG$ and $V_\ip:G\rightarrow \clG$ for $\ip=1,...,r$ are elements of the Lie algebra, denoted as $\clG$;  
$B_t^\ip$ and $W_t$ are mutually independent standard Wiener processes in $\Re$ 
and $\Re^m$, respectively, and they are also assumed to be independent of the initial state $X_0$;  
$h : G\rightarrow \Re^m$ is a given vector-valued nonlinear function.  
The $\io$-th coordinate of $Z_t$ and $h$ are denoted as $Z_t^\io$ and $h_\io$, 
respectively (i.e. $Z_t=(Z_t^1,...,Z_t^m)$ and $h=(h_1,...,h_m)$). 
The $\circ$ before $\ud B_t^\ip$ indicates that the stochastic differential equation (sde) 
\eqref{eq:dyn} is expressed in its Stratonovich form, which provides a coordinate-free description of the sde \cite{hsu2002}. 
The Einstein summation convention for the index $\ip$ is used in \eqref{eq:dyn}. 
A brief overview of matrix Lie groups and related notation is contained in \Sec{sec:prelim}.

\medskip

The problem is to numerically approximate the conditional distribution of $X_t$ given 
the time-history of observations $\UZ=\sigma(Z_s:s\leq t)$. 
The conditional distribution, denoted as $\pi_t^*$, acts on a function $f: G \rightarrow \Re$ 
according to 
\begin{equation*}
 \pi_t^*(f) := \Expect[f(X_t) | \UZ],
\end{equation*}
whose time-evolution is described by  
the Kushner-Stratonovich filtering equation (see Theorem 5.7 in \cite{xiong2008introduction}),  
\begin{align}
 \pi_t^*(f) = & ~\pi_0^*(f) + \int_0^t \pi_s^*(\Lstar f)\ud s + \int_0^t \big(\pi_s^*(fh)-\pi_s^*(h)\pi_s^*(f) \big)^T \big(\ud Z_s - \pi_s^*(h)\, \ud s \big) ,
 \label{eq:K-S}
\end{align}
for all $f\in \Cinfc(G)$ (smooth functions with compact support), 
where $\Lstar f := V_0\cdot f + \frac{1}{2} \sum_{\ip=1}^r V_\ip\cdot (V_\ip\cdot f)$.
The operations $V_0\cdot f$ and $V_\ip\cdot f$ is defined in \Sec{sec:prelim}.

\subsection{Literature Review}

Filtering of stochastic processes in non-Euclidean spaces has a rich history; c.f., 
\cite{duncan1977, caines1985IMA}. 
In recent years, the focus has been on computational approaches to numerically approximate the conditional distribution. 
Such approaches have been developed, e.g., by extending the classical extended Kalman filter (EKF) 
to Lie groups, and the extensions have appeared in discrete-time 
\cite{bourmaud2013discrete, zanetti2009norm}, continuous-time \cite{bonnabel2009IEKF, forbes2014continuous}, 
and continuous-discrete-time settings \cite{barrau2015TAC, bourmaud2015EKF}. 
In particular, a number of EKF-based filters have been proposed and applied for attitude estimation, 
e.g., the additive EKF \cite{bar1985attitude} and the multiplicative EKF \cite{lefferts1982kalman}. 
The EKF-based attitude filters require a linearized model of the estimation error, 
typically derived using one of the many 
three-dimensional attitude representations, e.g. the Euler angle \cite{itzhack87}, 
the rotation vector \cite{pittelkau2003rotation}, 
and the modified Rodrigues parameter \cite{markley2003attitude}. 

Apart from the EKF, particle filters for matrix Lie groups has been an
active area of research~\cite{chiuso2000, kwon2007particle,
  marjanovic2016engineer}.  Typically, particle filters adopt
discrete-time description of the dynamics and are based on importance
sampling and resampling numerical procedures.  For the attitude estimation problem, 
the unscented quaternion estimator \cite{crassidis2003unscented} 
and the bootstrap particle filter \cite{cheng2010particle,
  oshman2006attitude} have been developed, using one of the attitude
representations.  Other non-parametric approaches include filters based on 
certain variational formulations on the Lie groups \cite{zamani2013tac, izadi2014rigid}.

More recently, geometric group-theoretic methods for Lie groups have been developed. 
Deterministic nonlinear observers that respect the intrinsic geometry of the Lie groups appear in   
\cite{mahony2008TAC, lageman2010gradient, vasconcelos2010nonlinear}. 
A class of symmetry-preserving observers have been proposed to 
exploit certain invariance properties 
\cite{bonnabel2009observer}, leading to 
the invariant EKF \cite{bonnabel2009IEKF, barrau2015TAC}, 
the invariant unscented Kalman filter \cite{condomines2013IUKF}, 
the invariant ensemble Kalman filter \cite{barrau2015TAC}, 
and the invariant particle filter \cite{barrau2014cdc} algorithms within the stochastic filtering framework. 
A closely related theme is the use of 
non-commutative harmonic analysis for characterizing error propagation 
and Bayesian fusion on Lie groups \cite{chirikjian2016harmonic, wang2006error, wolfe2011bayesian}.
More comprehensive surveys can be found in~\cite{crassidis2007survey, zamani2013thesis}.


\subsection{Overview of the Paper}

The objective of this paper is to obtain a generalization of the feedback particle filter (FPF), 
originally developed in \cite{Tao_Thesis, Tao_Automatica, Tao_TAC} in the Euclidean settings, to the filtering problem \eqref{eq:dyn}-\eqref{eq:obs}. 
The main result is to show that the update formula in the original Euclidean setting carries over 
to the manifold setting. 

\medskip

The contributions of this paper are as follows: 

\smallskip

\noindent {$\bullet$ \bf Feedback particle filter on matrix Lie groups.} 
The extension of the FPF for matrix Lie groups is derived. The particle dynamics, expressed 
in their Stratonovich form, are shown to provide a coordinate-free 
description of the filter that automatically satisfies the geometric constraints of the manifold. 
Even in the manifold setting, the FPF is 
i) shown to admit an error-correction gain-feedback structure,  
and ii) proved to be an exact algorithm. Exactness means that, in the limit of large number 
of particles, the empirical distribution of the particles exactly matches the 
posterior distribution.

\smallskip

\noindent {$\bullet$ \bf Poisson equation on matrix Lie groups.}
The FPF algorithm requires numerical approximation of the gain function as a solution 
to a linear Poisson equation on the Lie group. 
An existence-uniqueness result for the solution is described in the Lie group setting. 
Two numerical methods are proposed to approximate the solution: 
i) In the Galerkin scheme, the gain function is approximated using a set of pre-defined basis functions; 
ii) In the kernel-based scheme, a numerical solution is obtained by 
solving a certain fixed-point equation.  

\medskip

\noindent {$\bullet$ \bf Feedback particle filter for attitude estimation.}
The attitude estimation problem represents the special case where the Lie group is $SO(3)$. 
For this important special case, 
the explicit form of the filter is described with respect to both the rotation matrix and the quaternion coordinates, 
with the latter being demonstrated for computational purposes. 
The FPF is also compared with some popular attitude filters, including 
the multiplicative EKF, the unscented quaternion estimator, the left invariant EKF, 
and the invariant ensemble Kalman filter. 
Simulation studies are presented to illustrate the performance comparison of these filters and the FPF algorithm.

\medskip

\noindent {$\bullet$ \bf The FPF with concentrated distributions.}  
The primary challenge in implementing the FPF algorithm arises due to the gain function approximation. 
In a certain special case, namely where the posterior distribution is concentrated, 
certain closed-form approximation, referred to as the {\em constant gain approximation}, of the gain function is obtained. 
For this approximation, evolution equations for the mean and the covariance are also derived and shown to 
be closely related to the left invariance EKF algorithm.

\subsection{Organization of the Paper}

The outline of the remainder of this paper is as follows:
A brief review  of the relevant Lie group preliminaries is contained in \Sec{sec:prelim}. 
In \Sec{sec:FPF_LG}, the generalization of the FPF algorithm to matrix Lie groups is presented, 
including both theory and numerical algorithms. 
The FPF algorithm for attitude estimation and its special case with concentrated conditional distributions are  
described in \Sec{sec:attitude_SO3} and \Sec{sec:FPF_CD}, respectively. 
Numerical simulations are provided in \Sec{sec:simulations}. 
The proofs appear in the Appendix.

\section{Mathematical Preliminaries}
\label{sec:prelim}

This section includes a brief review of matrix Lie groups.
The intent is to fix the notation used in subsequent sections.

\subsection{Geometry of Matrix Lie Groups}

The {\em general linear group}, denoted as $GL(n;\Re)$, is the group of $n\times n$ invertible matrices, 
where the group operations are the matrix multiplication and inversion. 
The identity element is the identity matrix, denoted as $I$.
A {\em matrix Lie group}, denoted as $G$, is a closed subgroup of $GL(n;\Re)$. 
$G$ is assumed to be connected. 
The {\em Lie algebra} of $G$, denoted as $\clG$, is the set of matrices $V$
such that the matrix exponential, $\exp (tV)$, is in $G$ for all $t \in \Re$. $\clG$ is a vector space whose 
dimension, denoted as $d$, equals the dimension of the group. 
The Lie algebra is equipped with an inner product, denoted as $\lr{\cdot}{\cdot}_\clG$, 
and an orthonormal basis $\{E_1,...,E_d\}$ with $\lr{E_i}{E_j}_\clG = \delta_{ij}$. 
The norm of $V \in \clG$ is defined as $\normG{V}:=\sqrt{\lr{V}{V}_\clG}$.

\medskip

{\em Example:} The special orthogonal group $SO(3)$ is the group of $3\times3$ 
matrices $R$ such that $RR^T=I$ and $\det(R)=1$. The Lie algebra $so(3)$ is the 3-dimensional 
vector space of skew-symmetric matrices. An inner product is 
$\lr{\Omega_1}{\Omega_2}_{so(3)} = \frac{1}{2}\,\Tr(\Omega_1^T\Omega_2)$ for $\Omega_1, \Omega_2 \in so(3)$, 
and an orthonormal basis $\{E_1,E_2,E_3\}$ of $so(3)$ is given by,
\begin{equation}
  \setlength{\arraycolsep}{2.9pt}
  E_1 = 
  \begin{bmatrix}
  0 & 0 & 0 \\
  0 & 0 & -1 \\
  0 & 1 & 0
  \end{bmatrix},~
  E_2 = 
  \begin{bmatrix}
  0 & 0 & 1 \\
  0 & 0 & 0 \\
  -1 & 0 & 0
  \end{bmatrix},~
  E_3 = 
  \begin{bmatrix}
  0 & -1 & 0 \\
  1 & 0 & 0 \\
  0 & 0 & 0
  \end{bmatrix}.
  \label{eq:basis_so3}
\end{equation}
These matrices have the physical interpretation of generating rotations about the 
three canonical axes in $\Re^3$. Here, $\det(\cdot)$ and $\Tr(\cdot)$ denote the 
determinant and trace of a matrix, respectively. Given the basis in \eqref{eq:basis_so3}, 
a vector $\omega = (\omega_1,\omega_2,\omega_3) \in \Re^3$ is uniquely mapped to an element 
in $so(3)$, denoted as $\Omega = \sk{\omega}:=\omega_1E_1 + \omega_2E_2 + \omega_3E_3$. 
Conversely, $\omega := \ve{\Omega}$.
\qed

\medskip

\begin{figure}
  \centering
    \includegraphics[width=.4\textwidth]{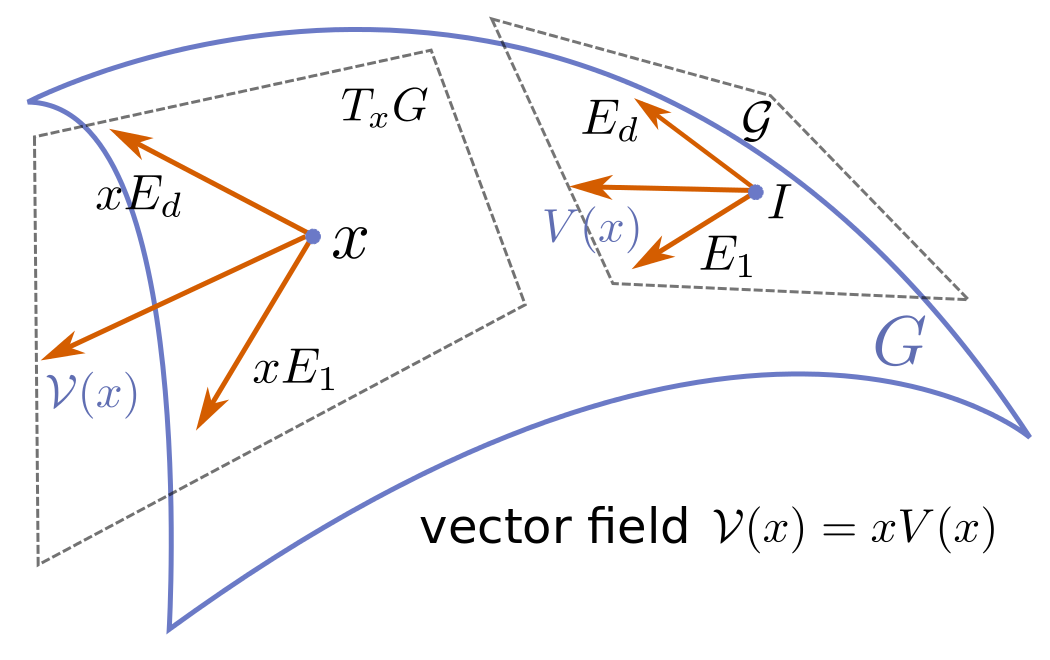}
  \caption{Construction of vector-fields on $G$.}
  \label{fig:vector_field}
\end{figure}

The Lie algebra is identified with the tangent space at the
identity matrix $I$, and can furthermore be used to construct a basis
$\{E_1^x,...,E_d^x\}$ for the tangent space at $\x\in G$, where $E_n^x
= \x\,E_n$ for $n=1,...,d$.  Therefore, a vector-field on $G$, denoted as $\clV$, is expressed as
  $$ \clV(\x) = v_1(\x) \, E_1^x + \cdots + v_d(\x) \, E_d^x, $$
with $v_n(x): G \rightarrow \Re$ for $n=1,...,d$. We write $\clV(x) = x\,V(x)$, 
where $V(x) := v_1(x)\,E_1 + \cdots + v_d(x)\,E_d$ is an element of the Lie algebra $\clG$ 
for each $\x\in G$. The functions $\big(v_1(\x), ..., v_d(\x)\big)$ are 
referred to as the {\em coordinates} of the vector-field. 
The construction of vector-fields on $G$ is illustrated in Figure \ref{fig:vector_field}.

The inner product $\lr{\cdot}{\cdot}_\clG$ in the Lie algebra induces an inner product of two vector-fields,
\begin{equation*}
 \lr{\clV}{\clW}(\x) := \lr{V(x)}{W(x)}_\clG = \sum_{n=1}^d v_n(\x)w_n(\x),
\end{equation*}
along with the associated norm $|\clV|_G(x):=\sqrt{\lr{\clV}{\clV}(x)}$.

With a slight abuse of notation, the action of the vector-field $\clV$ on a function $f: G \rightarrow \Re$ is denoted as,
\begin{equation}
 V\cdot f(\x) := \frac{\ud}{\ud \tau}\Big|_{\tau=0} f \big(\x \, \exp(\tau \, V(\x))\big).
 \label{eq:action}
\end{equation}

We also define the vector-field, $\grad(\phi)$, for a function $\phi: G \rightarrow \Re$ as,
\begin{equation}
 \grad(\phi)(\x) := \x \, \K(\x),
 \label{eq:gradient}
\end{equation}
where $\K(\x)\in \clG$ with coordinates $\big(E_1\cdot\phi(\x),...,E_d\cdot\phi(\x)\big)$ 
and according to \eqref{eq:action}, 
$(E_n\cdot \phi)(\x) := \frac{\ud}{\ud \tau} \big |_{\tau=0} \phi\big(x \,\exp(\tau \, E_n)\big)$ for $n=1,...,d$.
The vector-field acts on a function $f$ as,
\begin{equation}
 \K\cdot f(\x) = \lr{\grad(\phi)}{\grad(f)}(\x) = \sum_{n=1}^d E_n\cdot\phi(\x) \, E_n\cdot f(\x). 
 \label{eq:grad_inner_product}
\end{equation}

\medskip

We consider the following function spaces: 
The vector space of smooth real-valued functions $f:G\rightarrow\Re$ with compact support is denoted as $\Cinfc(G)$.
For a probability distribution $\pi$ on $G$, $L^2(G;\pi)$ denotes the Hilbert space of 
functions on $G$ that satisfy $\pi(|f|^2) < \infty$ (\,here $\pi(|f|^2) := \int_G |f|^2 \ud \pi(x)$\,);
$H^1(G;\pi)$ denotes the Hilbert space of functions $f$ such that $f$ and 
$E_n\cdot f$ (defined in the weak sense) are all in $L^2(G;\pi)$; 
and $H_0^1(G;\pi) := \{\phi\in H^1(G;\pi)\,|\,\pi(\phi)=0 \}$.

\subsection{Quaternions}

Quaternions provide a computationally efficient coordinate representation for $SO(3)$. 
A unit quaternion has the form,
\begin{align}
  q & = (q_0, ~q_1, ~q_2, ~q_3) = \Big( \cos(\frac{\theta}{2}), ~\omega_1\,\sin(\frac{\theta}{2}), ~\omega_2\,\sin(\frac{\theta}{2}), ~
  \omega_3\,\sin(\frac{\theta}{2}) \Big), 
  \label{eq:quaternion_def}
\end{align}
which represents rotation of angle $\theta$ about the axis defined 
by the unit vector $(\omega_1, \omega_2, \omega_3)$.
As with $SO(3)$, the space of quaternions admits a Lie group structure: 
The identity quaternion is denoted as $q_I:=(1,0,0,0)$, 
the inverse of $q$ is denoted as $q^{-1}:=(q_0,-q_1,-q_2,-q_3)$, and the multiplication is defined as,
\begin{equation*}
 p \otimes q := 
 \begin{bmatrix}
  p_0q_0 - p_V \cdot q_V \\
  p_0q_V + q_0p_V + p_V \times q_V
 \end{bmatrix},
\end{equation*}
where $p_V = (p_1,p_2,p_3)$, $q_V = (q_1,q_2,q_3)$, and $\cdot$ and $\times$ denote the 
dot product and the cross product between two vectors.

Given a unit quaternion $q$, the corresponding rotation matrix $R \in SO(3)$ is calculated by, 
\begin{equation}
 R = 
 \begin{bmatrix}
  2q_0^2 + 2q_1^2 - 1 & 2(q_1q_2-q_0q_3) & 2(q_1q_3+q_0q_2) \\
  2(q_1q_2+q_0q_3) & 2q_0^2 + 2q_2^2 - 1 & 2(q_2q_3-q_0q_1) \\
  2(q_1q_3-q_0q_2) & 2(q_2q_3+q_0q_1) & 2q_0^2 + 2q_3^2 - 1
 \end{bmatrix}.
\label{eq:convert_qR}
\end{equation}

\medskip

For more comprehensive introduction of matrix Lie groups and quaternions, we refer the reader to 
\cite{hall2015lie, trawny2005indirect}.

\section{Feedback Particle Filter on Matrix Lie Groups}
\label{sec:FPF_LG}

This section extends the FPF algorithm originally proposed in \cite{Tao_Automatica} to matrix
Lie groups, with necessary modifications to the original framework to account 
for the manifold structure.

\subsection{Particle Dynamics and Control Architecture}

The FPF on a matrix Lie group $G$ is a controlled system 
comprising of $N$ stochastic processes $\{X_t^i\}_{i=1}^N$ with $X_t^i \in G$. 
The particles are modeled by the Stratonovich sde,
\begin{align}
 \ud X_t^i = & ~X_t^i \big(V_0(X_t^i)+ \u(X_t^i,t)\big) \, \ud t + X_t^i \, V_\ip(X_t^i) \circ\ud B_{t}^{\ip,i} + X_t^i \, \K_\io(X_t^i,t) \circ \ud Z_t^\io,
 \label{eq:particle_dyn_uK}
\end{align}
where $B_t^{\ip,i}$ for $\ip = 1,...,r$ and $i=1,...,N$ are mutually independent standard Wiener processes in $\Re$, 
and the Einstein summation convention is used for the indices $\ip$ and $\io$. 
The functions $\u(x,t), ~\K_\io(x,t) :G \times [0,T] \rightarrow \clG$ are 
referred to as the {\em control} and {\em gain} function, respectively, 
whose coordinates are denoted as $(\u_1,...,\u_d)$ and $(\k_{1,\io},...,\k_{d,\io})$, for $j=1,...,m$. 
These functions need to be chosen. The following admissibility requirement is imposed on $\u$ and $\K_\io$:

\medskip

\begin{definition}
 {\em (Admissible Input)}:
 The functions $\u(x,t)$ and $\K_\io(x,t)$ are {\em admissible} if 
 they are $\UZ-$measurable and 
 $\Expect[\bpar{\sum_n|\u_n(X_t^i,t)|}]<\infty$,  
 $\Expect[\sum_{n} |\k_{n,\io}(X_t^i,t)|^2]<\infty$ for each $\io = 1,...,m$ and for all $t$. \qed
 \label{def:admissible}
\end{definition}

\medskip

The conditional distribution of the particle $X_t^i$ given $\UZ$ is denoted by $\pi_t$, 
which acts on a function $f$ according to
 $$ \pi_t(f) := \Expect[f(X_t^i) | \UZ]. $$
The evolution equation for $\pi_t$ is given by the proposition below. The proof appears in Appendix
\ref{apdx:prop_evolution}.

\medskip

\begin{proposition}
 Consider the particle $X_t^i$ with dynamics described by \eqref{eq:particle_dyn_uK}. 
 The forward evolution equation of the conditional distribution $\pi_t$ is given by,  
 \begin{equation}
  \pi_t(f) = \pi_0(f) 
   + \int_0^t \pi_s( \mathcal{L}f) \ud s 
   + \int_0^t \pi_s( \K_\io\cdot f) \ud Z_s^\io,
  \label{eq:forward}
 \end{equation}
for any $f\in\Cinfc(G)$, where the operator ${\cal L}$ is defined as
\begin{equation*}
 {\cal L}f := (V_0+\u)\cdot f + \frac{1}{2} \,\sum_{\ip=1}^r V_\ip\cdot(V_\ip\cdot f) 
 + \frac{1}{2} \sum_{\io=1}^m \K_\io\cdot (\K_\io\cdot f). 
\end{equation*}
\label{prop:evolution} 
\qed
\end{proposition}

\medskip

\noindent {\bf Problem statement:}
There are two types of conditional distributions of interest:
\begin{itemize}
 \item $\pi_t^{*}$: The conditional distribution of $X_t$ given $\UZ$.
 \item $\pi_t$: The conditional distribution of $X_t^i$ given $\UZ$.
\end{itemize}
The functions $\u(x,t),\,\K_\io(x,t)$ are said to be {\em exact} if $\pi_t=\pi_t^*$ 
for all $t \in [0,T]$. Thus, the objective is to choose $u$ and $\K_\io$ such that, 
given $\pi_0=\pi_0^{*}$, the evolution of the two conditional 
distributions are identical (see \eqref{eq:K-S} and \eqref{eq:forward}).

\medskip

\noindent{\bf Solution:}
The FPF represents the following choice of the gain function $\K$ and 
the control function $\u$:

\medskip

\noindent {\em 1) Gain function:} The gain function is obtained as follows: 
For $\io=1,...,m$, let $\phi_\io\in H^1(G;\pi_t)$ be the solution of a linear Poisson equation:
\begin{equation}
 \begin{aligned}
  & \pi_t \big( \lr{\grad(\phi_\io)}{\grad(\psi)} \big) = \pi_t \big( (h_\io-\hat{h}_\io) \psi \big), \\
  & \pi_t (\phi_\io) = 0 ~~~~(\text{normalization}), 
 \end{aligned}
 \label{eq:BVP}
\end{equation}
for all $\psi\in H^1(G;\pi_t)$, where $\hat{h}_\io := \pi_t(h_\io)$. The gain function $\K_\io \in \clG$ is then chosen as,
\begin{equation}
 \x \, \K_\io(\x, t) = \grad(\phi_\io)(\x). 
 \label{eq:gain_def}
\end{equation}
Given a basis $\{E_n\}_{n=1}^d$ of the Lie algebra $\clG$, 
and noting that (see \eqref{eq:gradient}) 
  $$ \grad(\phi_\io)(\x) = E_1\cdot\phi_\io(\x) \, E_1^x + \cdots + E_d\cdot\phi_\io(\x) \, E_d^x $$
where $E_n^x = x \, E_n$, the coordinates of $\K_\io$ is given by
\begin{equation}
 \k_{n,\io}(\x,t) = E_n\cdot \phi_\io(\x)~,~~\text{for}~ n=1,...,d.
 \label{eq:gain_coord}
\end{equation}

\medskip

\noindent {\em 2) Control function:} The function $\u$ is chosen as,
\begin{equation}
 \u(x,t) = -\frac{1}{2}\,\sum_{\io=1}^m \K_\io(x,t)\,\big(h_\io(x)+\hat{h}_\io\big).
 \label{eq:optimal_u}
\end{equation}

\noindent{\bf Feedback particle filter:}
Using these choice of $\u$ and $\K$, the $i$-th particle in the FPF has the following representation: 
\begin{equation}
 \ud X_t^i = \underbrace{X_t^i \, V_0(X_t^i) \, \ud t + \,X_t^i \, V_\ip \circ \ud B_{t}^{\ip,i}}_{\text{propagation}} 
 + \underbrace{X_t^i \,\, \K_\io(X_t^i,t) \circ \ud \innov_{t}^{\io,i}}_{\text{observation update}},
 \label{eq:particle_dyn}
\end{equation}
where the error $\ud\innov_t^{\io,i} \in \Re$ is a modified form of the innovation process:
\begin{equation}
 \ud \innov_{t}^{\io,i} = \ud Z_{t}^\io-\frac{h_\io(X_t^i)+\hat{h}_\io}{2}\ud t, 
 \label{eq:innovation}
\end{equation}
for each $\io=1,...,m$. 
The $i$-th particle implements the Bayesian update step\,--\,to account for the conditioning due to 
the observations\,--\,as gain times an error, which is akin to the 
feedback structure in a classical Kalman filter. 

Note that the Poisson equation \eqref{eq:BVP} must be solved for each $j=1,...,m$, 
and for each time $t$. 

The exactness is asserted in the following theorem. The proof 
is contained in appendix 
\ref{apdx:thm_consistency}.

\medskip

\begin{theorem}
 Let $\pi_t^{*}$ and $\pi_t$ satisfy the forward evolution equations 
 \eqref{eq:K-S} and \eqref{eq:forward}, respectively. 
 Suppose that the gain functions $\K_\io$, $\io=1,...,m$, are 
 obtained using \eqref{eq:BVP}-\eqref{eq:gain_def}, and the control function $\u$ is 
 obtained using \eqref{eq:optimal_u}.
 Suppose also that these functions are admissible. Then, assuming $\pi_0=\pi_0^{*}$, we have, 
 $$ \pi_t(f) = \pi_t^{*}(f),$$
 for all $t \in [0,T]$ and all function $f\in\Cinfc(G)$.
\label{thm:consistency} 
\qed
\end{theorem}

\medskip

\begin{remark}
In the original Euclidean setting, the FPF has the prettiest\,--\,gain
times error\,--\,representation of the update step in the Stratonovich
form of the filter (see Remark 1 in \cite{Tao_Automatica}).  In the \ito~form, the filter includes an
additional Wong-Zakai correction term. 
For sdes on a manifold, it is well known that the Stratonovich form is invariant
to coordinate transformations while the \ito~form is not \cite{oksendal2003}.  
So, for the gain times error form of the update step to have an intrinsic
coordinate independent form, the multiplication must necessarily be in
the Stratonovich form. The upshot is that the gain and error formula for the update step 
in \eqref{eq:particle_dyn} is intrinsic to the filter. 
\qed
\end{remark}

\medskip

\begin{remark}
The equation~\eqref{eq:BVP} is the weak form of a Poisson equation. 
Suppose $\pi_t$ admits an everywhere positive density, denoted as $\rho_t$.  Then the strong
form of~\eqref{eq:BVP} is given by the standard Poisson equation,
 \begin{equation}
  \lap \phi = -(h_\io-\hat{h}_\io),
\label{eq:PE_strong_form}
 \end{equation}
 where $\lap \phi := \frac{1}{\rho_t} \, \div \bpar{\rho_t \,
 \grad(\phi)}$ is the weighted Laplacian on the manifold \cite{grigoryan2009heat}, 
 and $\div(\cdot)$ denotes the divergence operator. 
 Multiplying both sides of~\eqref{eq:PE_strong_form} by $\psi(x)\rho_t(x)$ and integrating
 by parts, one arrives at the weak form~\eqref{eq:BVP}. 
 \qed
 \label{remark:gain}
 \end{remark}
 
\medskip
  
In the Euclidean case, the gain function was obtained as the
gradient of the solution of a Poisson equation~\cite{Tao_Automatica}. 
Remark \ref{remark:gain} shows that the Euclidean gain function is a
special case of the more general Lie group formula \eqref{eq:BVP}-\eqref{eq:gain_def}. 
For the latter, the definition of divergence and gradient 
ensures that the Poisson equation has a coordinate-free
representation. The gain function, expressed as gradient of the
solution of the Poisson equation, is an element of the Lie algebra.
This is consistent with the use of Lie algebra to define the vector
fields for dynamics evolving on the Lie group.   

In summary, the FPF is an intrinsic algorithm. 
The FPF update formula not only provides for a
generalization of the Kalman filter to the nonlinear non-Gaussian case
but also that the generalization carries over to nonlinear spaces such
as the Lie groups.  This is expected to have important consequences
for many applications in vision and robotics where Lie groups naturally arise.

\subsection{Well-posedness and Admissibility of the Gain}
\label{sec:gain}

The admissibility of the gain function solution, i.e., $\Expect[\bpar{\sum_n|\u_n(X_t^i,t)|}]<\infty$ and  
 $\Expect[\sum_{n} |\k_{n,\io}(X_t^i,t)|^2]<\infty$, requires a
well-posedness analysis of the Poisson equation.  As in the original
Euclidean setting, we make the following assumptions:

\medskip

\begin{assumption}
 The function $h_\io \in L^2(G;\pi_t)$ for each $\io = 1,...,m$ and for all $t$. 
 \qed
 \label{assumption:1}
\end{assumption}

\begin{assumption}
 The distribution $\pi_t$ admits a uniform spectral gap 
 (or Poincar\'{e} inequality) with constant $\bar{\lambda}$ (Sec. 4.2 in \cite{bakry2013analysis}): 
 That is, for a function $\phi \in H_0^1(G;\pi_t)$ 
 and for all $t\in[0,T]$,
 \begin{equation*}
  \pi_t\big( |\phi|^2 \big) \leq \frac{1}{\bar{\lambda}}\,\pi_t\big(|\grad(\phi)|^2 \big). 
 \end{equation*} 
 \qed
 \label{assumption:2}
\end{assumption}

The proof of the following well-posedness theorem is identical to the
proof presented in~\cite{Tao_Automatica} for the Euclidean case.  It is omitted. 

\medskip

\begin{theorem}
 Under Assumption \ref{assumption:1} and Assumption \ref{assumption:2}, the Poisson equation 
 \eqref{eq:BVP} possesses a unique solution $\phi_\io \in H_0^1(G;\pi_t)$, satisfying
 \begin{equation*}
  \pi_t \big( |\grad(\phi_\io)|_G^2 \big) \leq 
  \frac{1}{\bar{\lambda}} \, \pi_t \big( |h_j - \hat{h}_j|^2 \big),
 \end{equation*}  
 for each $j=1,...,m$ and for all $t$.  For this solution, one has the
 following bounds,
 \begin{align*}
  \pi_t\big( |\K_\io|_\clG^2 \big) & \leq \frac{1}{\bar{\lambda}}\,\pi_t\big( |h_j - \hat{h}_j|^2 \big),
  ~~~j = 1,...,m, \\
  \pi_t\big(\,\sum_{n=1}^d|\u_n|\,\big) & \leq C 
  \sum_{\io=1}^m \pi_t\big( |h_\io|^2 \big), 
 \end{align*}
 where the constant $C$ depends on $\bar{\lambda}$.  That is, the resulting gain and 
 control functions are admissible according to Definition \ref{def:admissible}.
 \qed
 \label{thm:E_U}
\end{theorem}

\medskip

\begin{remark}[Remark on Assumptions A1-A2]
Suppose the Lie group $G$ is compact, e.g., $SO(3)$.  In this case, if
$\pi_t$ has an everywhere positive density $\rho_t$, then Assumption \ref{assumption:1} and
\ref{assumption:2} automatically hold.  For non compact manifolds, e.g., $SE(3)$, the
assumptions hold if the density $\rho_t$ has a Gaussian tail (see
Remark 2 in~\cite{Tao_Automatica}). 
\qed
\end{remark}

\medskip
The main challenge to implement the FPF algorithm is to approximate
the gain function solution.  Since the problem~\eqref{eq:BVP}  is linear, the
approximation involves constructing a matrix problem to obtain the
approximate solution.    
In the following two sections, two numerical schemes for the approximation
are presented.  
Since the equations for each $\io=1,...,m$ are uncoupled, without loss of generality, 
a scalar-valued observation is assumed (i.e., $m=1$, and $\phi_\io$, $h_\io$, $\K_\io$ are 
denoted as $\phi$, $h$, $\K$).  
As the time $t$ is fixed, the explicit dependence on $t$ is suppressed 
(i.e., $\pi_t$, $\X_t^i$ are denoted as $\pi$, $\X^i$).

\subsection{Galerkin Gain Function Approximation}
\label{sec:galerkin_gain}

In a Galerkin approach, the solution $\phi$ is approximated as,
 $$ \phi = \sum_{l=1}^L \kappa_l \, \psi_l, $$
where $\{\psi_l\}_{l=1}^L$ is a given (assumed) set of {\em basis functions} on $G$. 
The coordinates of the gain function $\K$ with respect to a basis $\{E_n\}_{n=1}^d$ of $\clG$ are 
then given by (see \eqref{eq:gain_coord}),
  $$ \k_n = \sum_{l=1}^L \kappa_l \, E_n\cdot \psi_l, ~~~n = 1,...,d. $$

The finite-dimensional approximation of the Poisson equation \eqref{eq:BVP} 
is to choose coefficients $\{\kappa_l\}_{l=1}^L$ such that,
\begin{equation}
  \sum_{l=1}^L \kappa_l \, \pi\big( \lr{\grad(\psi_l)}{\grad(\psi)} \big)
     = \pi \big( (h-\hat{h})\psi \big),
  \label{eq:FEM}
\end{equation}
for all $\psi \in \text{span}\{\psi_1,...,\psi_L\} \subset H^1(G;\pi)$. 
On taking $\psi=\psi_1,...,\psi_L$, \eqref{eq:FEM} 
is compactly written as a linear matrix equation, 
\begin{equation} 
 A\kappa = b,
 \label{eq:kappa}
\end{equation}
where $\kappa := (\kappa_1,\hdots,\kappa_L)$, and  
the $L\times L$ matrix $A$ and the $L\times 1$ vector $b$ are defined and approximated as,
\begin{align*}
 A_{kl} & = \pi\big(\lr{\grad(\psi_l)}{\grad(\psi_k)} \big) \approx \frac{1}{N} \sum_{i=1}^N \lr{\grad(\psi_l)({\X}^i)}{\grad(\psi_k)({\X}^i)} \\
 & = \frac{1}{N} \sum_{i=1}^N \sum_{n=1}^d (E_n\cdot\psi_l)({\X}^i) \, (E_n\cdot\psi_k)({\X}^i), \\
 b_k & = \pi \big((h-\hat{h})\psi_k \big) 
 \approx \frac{1}{N}\sum_{i=1}^N (h({\X}^i)-\hat{h})\psi_k({\X}^i),
\end{align*}
where $\hat{h}\approx\frac{1}{N}\sum_{i=1}^N h({\X}^i) =: \hat{h}^{(N)}$.

\medskip

Note that both the Poisson equation \eqref{eq:BVP} as well as its Galerkin finite-dimensional approximation 
\eqref{eq:kappa} are coordinate-free representations. Particle-based
approximation of the solution 
\eqref{eq:kappa} can be carried out for any choice of coordinates.
Certain coordinates may offer computational advantages, e.g.,
quaternions for SO(3) (see Sec.~\ref{sec:attitude_SO3}).

The non-trivial step in the Galerkin approximation is the choice of
the basis function.  In general, this choice is problem dependent.
For matrix Lie groups, one choice is to use the Fourier basis.  

\medskip

{\em Basis functions on $SO(2)$:} 
The Lie group $SO(2)$ is identified with the unit circle $S^1$. Using the angle coordinate $\theta \in S^1$, 
the simplest choice of the basis functions are the Fourier basis, e.g., 
\begin{equation}
 \psi_1(\theta) = \sin(\theta) , ~~\psi_2(\theta) = \cos(\theta). 
 \label{eq:basis_SO2}
\end{equation}
Note that these two basis functions are the eigenfunctions of the Laplacian on $SO(2)$, 
associated with its smallest non-zero eigenvalue. 
\qed
  
\medskip

For more general Lie groups, e.g., $SO(3)$, the Fourier basis are the
eigenfunctions of the Laplace-Beltrami operator on the manifold \cite{lablee2015spectral}.
Given its importance in applications, the eigenfunctions associated
with the smallest eigenvalue for $SO(3)$ appears in Appendix \ref{apdx:basis_SO3}. 
Also included are the necessary calculations to implement the Galerkin procedure in this case.  

\medskip

The gain function approximation is a hard problem.  The Galerkin
algorithm represents the most straightforward solution where the
hard part is to select the basis functions. 
A number of papers have considered related approaches: 
i) the use of proper orthogonal decomposition (POD) to select basis functions in~\cite{berntorp2016acc}; 
ii) a continuation scheme in~\cite{matsuura2016suboptimal}; and 
iii) certain probabilistic approaches involving dynamic programming in \cite{Sean_CDC2016}.  We
expect that many of these approaches will also generalize to the manifold setting.  

In the following, we present a recent approach from~\cite{AmirCDC2016} 
whose attractive feature is that it does not involve selection of basis functions.
In the numerical results presented in \Sec{sec:simulations}, we show that this
approach is also very effective.

\subsection{Kernel-based Gain Function Approximation}
\label{sec:kernel_gain}

In~\cite{AmirCDC2016}, the unknown function $\phi(x)$\,--\,solution of the Poisson equation \eqref{eq:BVP}
\,--\,is approximated by its values at the particles $\{X^i\}_{i=1}^N$:
\[
\varPhi := \bpar{\phi(X^1),\phi(X^2),\hdots,\phi(X^N)}.
\]
In terms of $\varPhi$, the BVP~\eqref{eq:BVP} is approximated as a
finite-dimensional fixed-point problem,
\begin{equation}
 \varPhi = \Ten \varPhi + \epsilon H^{(N)},
 \label{eq:fixed_point_approx}
\end{equation}
on the co-dimension $1$ subspace of normalized (i.e., mean zero) vectors, where $\epsilon$ is a small positive parameter, 
$H^{(N)} := \big( h(\X^1)-\hat{h}^{(N)}, h(\X^2)-\hat{h}^{(N)},
...,h(\X^N)-\hat{h}^{(N)} \big) \in \Re^N$,  and $\Ten\in\Re^{N\times N}$ is a
Markov matrix that is assembled from the ensemble
$\{X^i\}_{i=1}^N$. It is shown in~\cite{AmirCDC2016} that: 
\begin{enumerate}
\item The Markov matrix $\Ten$ is a strict contraction on the subspace, and thus    
\item the finite-dimensional
problem~\eqref{eq:fixed_point_approx} admits a unique 
normalized solution $\varPhi$,
\item this solution can be obtained by successive approximations, and 
\item $\varPhi$ approximates the true solution $\phi$ as $\epsilon\rightarrow 0$ and $N \rightarrow \infty$. 
\end{enumerate} 

For the manifold, the $(i,j)^{\text{th}}$ element of the $N\times N$ matrix $\Ten$ is constructed as,
\begin{equation}
 \Ten_{ij} = \frac{\wken (\X^i, \X^j)}{\sum_{l=1}^N \wken(\X^i,\X^l)},
 \label{eq:Tij}
\end{equation}
where the kernel $\wken: G \times G \rightarrow \Re$ is given by,
\begin{equation}
 \wken (\X^i, \X^j) = \frac{\ke(\X^i, \X^j)}
 {\sqrt{\frac{1}{N}\sum_{l=1}^N \ke(\X^i, \X^l)} 
 \sqrt{\frac{1}{N}\sum_{l=1}^N \ke(\X^j, \X^l)}}, 
 \label{eq:modified_kernel}
\end{equation}
and $\ke$ is the Gaussian kernel, 
\begin{equation}
 \ke (\X^i, \X^j) := \frac{1}{(4\pi\epsilon)^{d/2}} \exp \Big( -\frac{\zeta^2(\X^i, \X^j)}{4\epsilon} \Big),
 \label{eq:kernel} 
\end{equation}
where $d$ is the dimension of $G$, 
and $\zeta : G \times G \rightarrow \Re$ denotes a distance metric on $G$ 
induced from the Euclidean space in which $G$ is smoothly embedded (see Assumption 19 in \cite{hein2006graph}).  

\medskip

\begin{remark}
 The justification of the fixed-point problem
 \eqref{eq:fixed_point_approx} is as follows: 
Any solution of the Poisson equation \eqref{eq:BVP} is equivalently also the solution of the fixed-point problem, 
 \begin{equation}
  \phi = e^{\,\epsilon\,\Delta_\rho}\phi + \int_0^\epsilon e^{\,s\,\Delta_\rho} (h-\hat{h})\,\ud s, 
  \label{eq:fixed_point}
 \end{equation}
 where $e^{\,\epsilon\,\Delta_\rho}$ denotes the subgroup of $\Delta_\rho$. 
 In the limit as $\epsilon \rightarrow 0$ and $N \rightarrow \infty$, $\Ten$ represents a finite-dimensional approximation of 
 $e^{\,\epsilon\,\Delta_\rho}$ (see Proposition 3 in \cite{coifman2006diffusion}).
 \qed
\end{remark}

\medskip

The coordinates of the gain function $\k_n = E_n\cdot\phi$ for $n=1,...,d$ 
are obtained by taking an explicit derivative of \eqref{eq:fixed_point}. 
The calculations are summarized below: 
\begin{enumerate}
 \item Define the vector 
   $$ \widetilde{H}_n := \bpar{E_n\cdot h(X^1),\,E_n\cdot h(X^2),...,\,E_n\cdot h(X^N)}, $$
 and define the $N \times N$ matrix $\widetilde{Z}_n$ whose elements are given by,
   $$ (\widetilde{Z}_n)_{ij} := E_n\cdot \zeta^2(X^i,X^j), $$
 where $E_n\cdot\zeta^2(x,y) := \frac{\ud}{\ud \tau}\big|_{\tau=0} \zeta^2\bpar{x\,\exp(\tau E_n), y}$ for $x,y \in G$. 
 \item Define the $N \times N$ matrix, 
   $$ S_n := \Ten * \widetilde{Z}_n, $$
 where $*$ denotes the Hadamard (element-wise) product of two matrices, i.e., $(S_n)_{ij} = (\Ten)_{ij}\,(\widetilde{Z}_n)_{ij}$. 
 \item Define $\varUpsilon_n:=\bpar{\k_n(X^1),\,\k_n(X^2),...,\,\k_n(X^N)} \in \Re^N$. Then,  
 \begin{equation}
  \varUpsilon_n = \epsilon\,\widetilde{H}_n - \frac{1}{4\,\epsilon} 
  \big[ S_n \varPhi - (S_n \mathbf{1}) * (\Ten\varPhi) \big], 
  \label{eq:kernel_gain_coord}
 \end{equation}
 where $\mathbf{1} = (1,1,...,1)\in\Re^N$, and $*$ denotes the Hadamard 
 product of two vectors. 
\end{enumerate}

\medskip

{\em Example:} On the Lie group $SO(3)$, $d=3$, and the distance metric is given by (see \cite{huynh2009metrics}), 
  $$ \zeta^2(R_1, R_2) = |R_1-R_2|_F^2, $$
for $R_1, R_2 \in SO(3)$, where $|\cdot|_F$ is the Frobenius norm of a matrix. 
This metric is induced from the Euclidean space $\Re^9$, where the smooth embedding $i: SO(3) \rightarrow \Re^9$ is 
defined as $i(R) = (R_{11},R_{12},...,R_{33})$. 
Using the basis of $so(3)$ given by \eqref{eq:basis_so3}, we have 
  $$ E_n\cdot(\zeta^2)(R^i,R^j) = -2\,\Tr(R^i E_n R^j) $$ 
for $n=1,2,3$. 
\qed

\subsection{FPF Algorithm Summary}
\label{sec:FPF_summary}

\begin{algorithm}[t]       
        \caption{Feedback Particle Filter on a matrix Lie group}
        \begin{algorithmic}[1]
            \STATE {\bf initialization:}  sample $\{\X_0^i\}_{i=1}^N$ from $\pi_0^{*}$
            \STATE Assign $t=0$
            \STATE {\bf iteration:} from $t$ to $t+\Dt$
            \STATE Calc. $\hat{h}_j^{(N)} = \frac{1}{N} \, \sum_{i=1}^N h_j(\X_{t}^{i})$ for $j=1,2,...,m$
	    \FOR{$i=1$ to $N$}
	      \STATE Generate samples $\Delta B_{t}^{\ip,i}$ from $N(0,\Delta t)$ for $\ip=1,...,r$
	      \STATE Assign $\Delta U_t^i = 0$
	      \FOR{$j=1$ to $m$}
		  \STATE Calc. the error 
		         $\Delta \innov_{\io,t}^i := \Delta Z_{\io,t} - \frac{1}{2} \big( h_\io(\X_t^i) + \hat{h}_\io^{(N)} \big) \, \Dt $
		  \STATE Calc. gain function $\K_\io(\X_t^i,t)$
		  \STATE Assign $\Delta U_t^i = \Delta U_t^i + \K_\io(\X_t^i,t)\,\Delta \innov_{\io,t}^i$
	      \ENDFOR
	      \STATE Calc. $ \Delta V_t^i = V_0(X_t^i) \, \Dt + V_\ip(X_t^i) \, \Delta B_t^{\ip,i}
	                         + \Delta U_t^i$
	      \STATE Propagate the particle  
		      $\X_{t+\Dt}^i = \X_t^i \, \exp\big( \Delta V_t^i \big)$ 
	    \ENDFOR
            \STATE {\bf return:} empirical mean of $\{\X_{t+\Dt}^i\}_{i=1}^N$ 
            \STATE Assign $t = t+\Dt$
        \end{algorithmic}
        \label{alg:FPF-LG}
\end{algorithm}

The numerical algorithm of the FPF on a matrix Lie group is summarized in Algorithm \ref{alg:FPF-LG}. 
The algorithm simulates $N$ particles, $\{\X_t^i\}_{i=1}^N$, 
according to the sde \eqref{eq:particle_dyn}  
with the initial conditions $\{\X_0^i\}_{i=1}^N$ sampled i.i.d. from a given prior distribution $\pi_0^{*}$. 
The gain function is approximated using either the Galerkin scheme (see \Sec{sec:galerkin_gain} 
and Algorithm \ref{alg:galerkin}) or the kernel-based scheme (see \Sec{sec:kernel_gain} and Algorithm \ref{alg:kernel}). 

\begin{algorithm}[H]        
        \caption{Galerkin gain function approximation}
        \begin{algorithmic}[1]
        \STATE {\bf input:} Particles $\{\X^i\}_{i=1}^N$, basis functions $\{\psi_l\}_{l=1}^L$
        \STATE Calc. $\hat{h}^{(N)} = \frac{1}{N} \sum_{i=1}^N h(\X^{i})$
        \FOR{$k=1$ to $L$}
            \STATE Calc. $b_k = \frac{1}{N} \sum_{i=1}^N \big( h(\X^i)-\hat{h}^{(N)} \big)\psi_k(\X^i)$
            \FOR{$l=1$ to $L$}
                \STATE Calc. $A_{kl} = 
                \frac{1}{N} \sum_{i=1}^N \sum_{n=1}^d (E_n\cdot\psi_l)(\X^i) \, (E_n\cdot\psi_k)(\X^i)$
            \ENDFOR
        \ENDFOR
        \STATE Solve the matrix equation $A\kappa=b$, with $A=[A_{kl}]$, $b=[b_k]$
        \STATE Calc. $\k_n(\X^i) = \sum_{l=1}^L \kappa_l \, E_n\cdot\psi_l(\X^i)$, for $n=1,...,d$
        \STATE {\bf return:} Coordinates $\big\{\big( \k_n(\X^1),\,...,\k_n(\X^N) \big)\big\}_{n=1}^d$
        \end{algorithmic}
        \label{alg:galerkin}
\end{algorithm}

\begin{algorithm}[H]        
        \caption{Kernel-based gain function approximation}
        \begin{algorithmic}[1]
        \STATE {\bf input:} Particles $\{\X^i\}_{i=1}^N$, parameters $\epsilon$, $K$
        \STATE Calc. $\hat{h}^{(N)} = \frac{1}{N} \sum_{i=1}^N h(\X^{i})$ 
        \STATE Calc. $H^{(N)}_i = h(\X^i) - \hat{h}^{(N)}$ for $i=1,...,N$
        \STATE Calc. $\ke(\X^i,\X^j)$, $\wken(\X^i,\X^j)$ by \eqref{eq:kernel}, \eqref{eq:modified_kernel} for all $i,j$
        \STATE Calc. $\Ten_{ij}$ according to \eqref{eq:Tij} for all $i,j$
        \STATE Assign initial condition $\varPhi_0$ 
        \FOR{$k=0$ to $K-1$}
	    \STATE Calc. $\varPhi_{k+1} = \Ten \varPhi_k + \epsilon H^{(N)}$, with $\Ten=[\Ten_{ij}]$
	    \STATE Assign $\varPhi_{k+1} = \varPhi_{k+1} - \frac{1}{N}\,\sum_{i=1}^N (\varPhi_{k+1})_i$
	\ENDFOR
	    \STATE Calc. $\varUpsilon_n = \big( \k_n(\X^1),\,...,\k_n(\X^N) \big)$ for $n=1,...,d$ 
	    according to \eqref{eq:kernel_gain_coord} with $\varPhi = \varPhi_K$
        \STATE {\bf return:} Coordinates $\big\{\big( \k_n(\X^1),\,...,\k_n(\X^N) \big)\big\}_{n=1}^d$
        \end{algorithmic}
        \label{alg:kernel}
\end{algorithm}

\section{Attitude Estimation with Feedback Particle Filter}
\label{sec:attitude_SO3}

This section considers the problem of attitude estimation, cast as a continuous-time 
nonlinear filtering problem on the Lie group $SO(3)$. The explicit form of the FPF algorithm is 
described with respect to both the rotation matrix and the quaternion coordinates, with the latter 
being demonstrated for computational purposes. 

\subsection{Problem Formulation}
\label{sec:attitude_problem}

\noindent{\bf Process model:} A kinematic model of rigid body is given by,
\begin{equation}
 \ud R_t = R_t \, \Omega_t \ud t + R_t \circ \sk{\sigma_B \, \ud B_t},
 \label{eq:kinematics}
\end{equation}
where $R_t \in SO(3)$ is the orientation of the rigid body at time $t$, 
expressed with respect to an inertial frame; 
$\Omega_t = \sk{\omega_t}$ where $\omega_t \in \Re^3$ represents the angular velocity expressed in 
the body frame; $B_t$ is a standard Wiener process in $\Re^3$, and 
$\sigma_B$ is a positive scalar. 
Both $\Omega_t$ and $\sk{\sigma_B \, \ud B_t}$ are elements of the Lie algebra $so(3)$. 

Using the quaternion coordinates, \eqref{eq:kinematics} is written as,
\begin{equation}
 \ud q_t = \frac{1}{2} \, q_t \otimes (\omega_t \ud t + \sigma_B \, \ud B_t), 
 \label{eq:kinematics_quat}
\end{equation}
where, by a slight abuse of notation, $\omega_t \in \Re^3$ is interpreted as a quaternion 
$(0,\omega_t)$, and $\ud B_t$ is interpreted similarly. The sde \eqref{eq:kinematics_quat} is 
also interpreted in the Stratonovich sense.

\medskip

\noindent {\bf Accelerometer:} 
In the absence of translational motion, the accelerometer is modeled as (see \cite{mahony2008TAC}), 
\begin{equation}
 \ud Z_t^{g} = -R_t^T r^g \ud t + \sigma_W \ud W_t^g,
 \label{eq:accelerometer}
\end{equation}
where $r^g \in \Re^3$ is the unit vector in the inertial frame aligned with the gravity, 
$W_t^g$ is a standard Wiener process in $\Re^3$, and a parameter $\sigma_W$ is used to scale the observation noise. 

\medskip

\noindent {\bf Magnetometer:}
The model of the magnetometer is of a similar form (see \cite{mahony2008TAC}),
\begin{equation}
 \ud Z_t^{b} = R_t^T r^b \ud t + \sigma_W \ud W_t^b,
 \label{eq:magnetometer} 
\end{equation}
where $r^b \in \Re^3$ is the unit vector in the inertial frame aligned with the local magnetic field, 
and $W_t^b$ is a standard Wiener process in $\Re^3$. 

\medskip

In terms of the process and observation models \eqref{eq:kinematics}-\eqref{eq:magnetometer}, 
the nonlinear filtering problem for 
attitude estimation is succinctly expressed as,
\begin{subequations}
 \begin{align}
  \ud R_t & = R_t \, \Omega_t \ud t + R_t \circ \sk{\sigma_B \, \ud B_t}, \label{eq:kinematics_1} \\
  \ud Z_t & = h(R_t) \ud t + \sigma_W \ud W_t, \label{eq:observation_1}
 \end{align}
\end{subequations}
where $h:SO(3) \rightarrow \Re^6$ is a given function whose $j$-th coordinate 
is denoted as $h_j$, and $W_t$ is a standard Wiener process in $\Re^6$. 
Note that \eqref{eq:observation_1} encapsulates the sensor models given in 
\eqref{eq:accelerometer} and \eqref{eq:magnetometer} within a single equation. 
It is assumed that $B_t$ and $W_t$ are mutually independent, and both are independent of the initial 
condition $R_0$.

\medskip

\begin{remark}
 There are a number of simplifying assumptions implicit in the model defined in 
 \eqref{eq:kinematics_1}-\eqref{eq:observation_1}. In practice, $\omega_t$ needs to be 
 estimated from noisy gyroscope measurements and there is translational motion as well. 
 This requires additional models which are easily incorporated within the proposed filtering framework.  
 The purpose here is to elucidate the geometric aspects 
 of the FPF in the simplest possible setting of $SO(3)$. More practical FPF-based filters 
 that also incorporate models for translational motion, measurements of $\omega_t$ 
 from gyroscope, effects of translational motion on accelerometer, and effects of sensor bias 
 are subject of separate publication.
 \qed
 \label{remark:gyro}
\end{remark}

\subsection{FPF for Attitude Estimation}
\label{sec:FPF_attitude}

Following the general framework of FPF , 
the dynamics of the $i$-th particle is defined by,
\begin{equation}
 \ud R_t^i = R_t^i \, \Omega_t \ud t + R_t^i \circ \sk{\sigma_B \, \ud B_t^i} + R_t^i \, \sk{\K(R_t^i,t) \circ \ud \innov_t^i},
 \label{eq:particle_dyn_SO3}
\end{equation}
where $B_t^i$ for $i=1,...,N$ are mutually independent standard Wiener processes in $\Re^3$. 
The error $\ud \innov_t^i \in \Re^6$ is given by, 
\begin{equation*}
 \ud \innov_t^i = \ud Z_t - \frac{1}{2} \big( h(R_t^i) + \hat{h} \big) \ud t.
\end{equation*}

The gain function $\K$ is a $3 \times 6$ matrix whose entries are obtained as follows: 
For $j = 1,2,...,6$, the $j$-th column of $\K$ contains the coordinates of the vector-field $\grad (\phi_j)$, 
where the function $\phi_j \in H^1(SO(3);\pi_t)$ is a solution to the Poisson equation,
\begin{equation}
 \begin{aligned}
  & \pi_t \big( \lr{\grad(\phi_j)}{\grad(\psi)} \big) = \frac{1}{\sigma_W^2} \, \pi_t \big( (h_j-\hat{h}_j) \psi \big), \\
  & \pi_t (\phi_j) = 0 ~~~~(\text{normalization}), 
 \end{aligned}
 \label{eq:BVP_SO3}
\end{equation}
for all $\psi \in H^1(SO(3);\pi)$. 

\medskip

For numerical purposes, it is convenient to express the FPF with respect to 
the quaternion coordinates. In this coordinate representation, the dynamics of the $i$-th particle is given by,
\begin{equation}
 \ud q_t^i = \frac{1}{2} \, q_t^i \otimes \ud \nu_t^i, 
 \label{eq:FPF_quat}
\end{equation}
where $q_t^i$ is the quaternion state of the $i$-th particle, and $\nu_t^i \in \Re^3$ evolves according to,
\begin{equation}
 ~~\ud\nu_t^i = \omega_t \ud t + \ud B_t^i + \K(q_t^i,t) 
 \circ \Big( \ud Z_t-\frac{h(q_t^i) + \hat{h}}{2}\, \ud t \Big),
 \label{eq:FPF_domega}
\end{equation}
where $\K(q,t)=\K(R(q),t)$ and $h(q)=h(R(q))$, with $R=R(q)$ given by the formula \eqref{eq:convert_qR}.

\subsection{FPF Algorithm Summary}
\label{sec:FPF_summary_attitude}

\begin{algorithm}[t]        
        \caption{Feedback Particle Filter for attitude estimation}
        \begin{algorithmic}[1]
            \STATE {\bf initialization:}  sample $\{q_0^i\}_{i=1}^N$ from $\pi_0^{*}$
            \STATE Assign $t=0$
            \STATE {\bf iteration:} from $t$ to $t+\Dt$
            \STATE Calc. $\hat{h}^{(N)} = \frac{1}{N} \sum_{i=1}^N h(q_{t}^{i})$
	    \FOR{$i=1$ to $N$}
	      \STATE Generate a sample, $\Delta B_t^{i}$, from $N\bpar{0,(\Dt)\id}$
	      \STATE Calc. the error
		     $\Delta \innov_{t}^i := \Delta Z_{t} - \frac{1}{2} \, \big( h(q_t^i) + \hat{h}^{(N)} \big) \, \Dt $
	      \STATE Calc. gain function $\K(q_t^i,t)$ using Galerkin or kernel-based scheme
	      \STATE Calc. $ \Delta \nu_t^i = \omega_t \, \Dt +  \sigma_B \, \Delta B_t^i + \K(q_t^i,t) \, \Delta \innov_{t}^i $
	      \STATE Propagate the particle $q_t^i$ according to (see \cite{trawny2005indirect}, and $|\cdot|$ denotes the Euclidean norm in $\Re^3$)
	             \begin{equation*}
		      q_{t+\Dt}^i = q_t^i \otimes 
		      \begin{bmatrix}
		      \cos \big( |\Delta\nu_t^i|/2 \big) \\
		      \frac{\Delta\nu_t^i}{|\Delta\nu_t^i|}\, \sin \big( |\Delta\nu_t^i|/2 \big)
		      \end{bmatrix}
		     \end{equation*}	             
	    \ENDFOR
	    \STATE Define matrix $Q = \frac{1}{N} \sum\nolimits_{i=1}^N q_{t+\Dt}^i {q_{t+\Dt}^{i~T}}$
            \STATE {\bf return:} empirical mean of $\{q_{t+\Dt}^i\}_{i=1}^N$, 
                                 i.e., the unit eigenvector of $Q$ associated with its largest eigenvalue 
            \STATE Assign $t = t+\Dt$
        \end{algorithmic}
        \label{alg:FPF_SO3}
\end{algorithm}

The FPF algorithm is numerically implemented using the quaternion coordinates, 
and is described in Algorithm \ref{alg:FPF_SO3}. 
The algorithm simulates $N$ particles, $\{q_t^i\}_{i=1}^N$, according to the sde's \eqref{eq:FPF_quat} and \eqref{eq:FPF_domega}, 
with the initial conditions $\{q_0^i\}_{i=1}^N$ sampled i.i.d. from a given prior distribution $\pi_0^{*}$. 
The gain function is approximated using either the Galerkin scheme (see Sec. \ref{sec:galerkin_gain}) 
with the basis functions given in Appendix \ref{apdx:basis_SO3}, 
or the kernel-based scheme (see Sec. \ref{sec:kernel_gain}). 

Given a particle set $\{q_t^i\}_{i=1}^N$, its empirical mean is obtained as the 
eigenvector (with norm 1) of the $4\times4$ matrix $Q = \frac{1}{N} \sum\nolimits_{i=1}^N q_t^i {q_t^{i~T}}$ 
corresponding to its largest eigenvalue 
\cite{markley2007averaging}.

\section{Feedback Particle Filter with Concentrated Distributions}
\label{sec:FPF_CD}

In its original Euclidean setting \cite{Tao_Automatica}, the FPF algorithm is shown to represent a generalization of the Kalman filter 
in the following sense: Suppose that the signal and the observation models are linear and that 
the prior distribution is Gaussian. Then, it is shown that: 
\begin{enumerate}
 \item The gain $\K_t$ is a constant for each $t$ whose value equals the Kalman gain; 
 \item The conditional distribution $\pi_t$ of $\X_t^i$ is Gaussian whose mean and covariance evolve 
according to the Kalman filter.
\end{enumerate}

For the general nonlinear non-Gaussian case, 
the gain function $\K_t$ is no longer a constant and must be numerically approximated. 
However, the conditional expectation of the gain function, $\Expect[\K_t\,|\,\UZ]$, admits a closed-form expression 
which can furthermore be approximated using only the particles. 
The resulting approximation is referred to as the constant gain approximation.  
This approximation reduces to the Kalman gain in the linear Gaussian case. For the general case,
this approximation often suffices in practice particularly so when the conditional distribution is
unimodal \cite{Tao_Automatica, berntorp2015, stano2014}. 

On a Riemannian manifold, unfortunately, even the state space does not possess a linear structure. 
However, under the additional assumption that the posterior distribution is {``concentrated''} (see \cite{wang2006error}), 
one can expect the results to be close to the Euclidean case. 
In this section, the following is shown for the special case of concentrated distributions on matrix Lie groups:
\begin{enumerate}
 \item A closed-form formula for the constant gain approximation is derived and shown to equal the Kalman gain; 
 \item The equation for the mean and covariance are derived and shown to be closely related to the continuous-time
left invariant EKF algorithm in \cite{bonnabel2009IEKF}.
\end{enumerate} 

In this section, we restrict our attention to the filtering problem 
\eqref{eq:kinematics_1}-\eqref{eq:observation_1} on $SO(3)$. 
Such a restriction is not necessary but leads to a simpler presentation without undue notational burden.
Also, it allows us to make comparisons with the literature on filters for attitude estimation.

\subsection{Constant Gain Approximation of FPF}
\label{sec:CGA}

Consider {\em concentrated distribution} whereby the random variable $R$ on $SO(3)$ 
is parametrized as,
  $$ R = \mu\,\exp(\epsilon\,\sk{\dev}), $$
where $\chi \in so(3) \cong \Re^3$ is a Gaussian random variable with mean $0$ and covariance $\cov$, 
and $\epsilon$ is a small parameter. Formally, most of the probability mass of a concentrated distribution is supported 
in a small neighborhood of $\mu$, and the analysis pertains to the consideration of the asymptotic limit as $\epsilon \rightarrow 0$. 

The following proposition provides an approximate formula for the gain in this special case. 
The proof appears in Appendix \ref{apdx:constant_gain}. 
For notational ease, the dependence on the time $t$ is suppressed (i.e., we express $R_t$ as $R$, $\dev_t$ as $\dev$, $\pi_t$ as $\pi$ etc.).

\medskip

\begin{proposition}
 Consider the Poisson equation \eqref{eq:BVP_SO3} where the random variable $R = \mu\,\exp(\epsilon\,\sk{\dev})$, 
 and $\dev \in \Re^3$ is a Gaussian random variable with mean $0$ and covariance $\cov$. 
 Suppose $\sigma_W = \epsilon \, \sigmaW$. 
 Let $[\K]_\io: SO(3) \rightarrow \Re^3$ denote the $\io$-th column of the gain function $\K$. 
 Then, in the asymptotic limit as $\epsilon \rightarrow 0$, 
 \begin{equation*}
  [\K]_\io = \frac{1}{\overline{\sigma}_W^2} \, \cov \, H_\io + O(\epsilon), 
 \end{equation*}
 where $H_\io := \bpar{ E_1\cdot h_\io(\mu),\,E_2\cdot h_\io(\mu),\,E_3\cdot h_\io(\mu) } \in \Re^3$.
 \qed
 \label{prop:constant_gain}
\end{proposition}

\begin{figure}
  \centering
    \includegraphics[width=.44\textwidth]{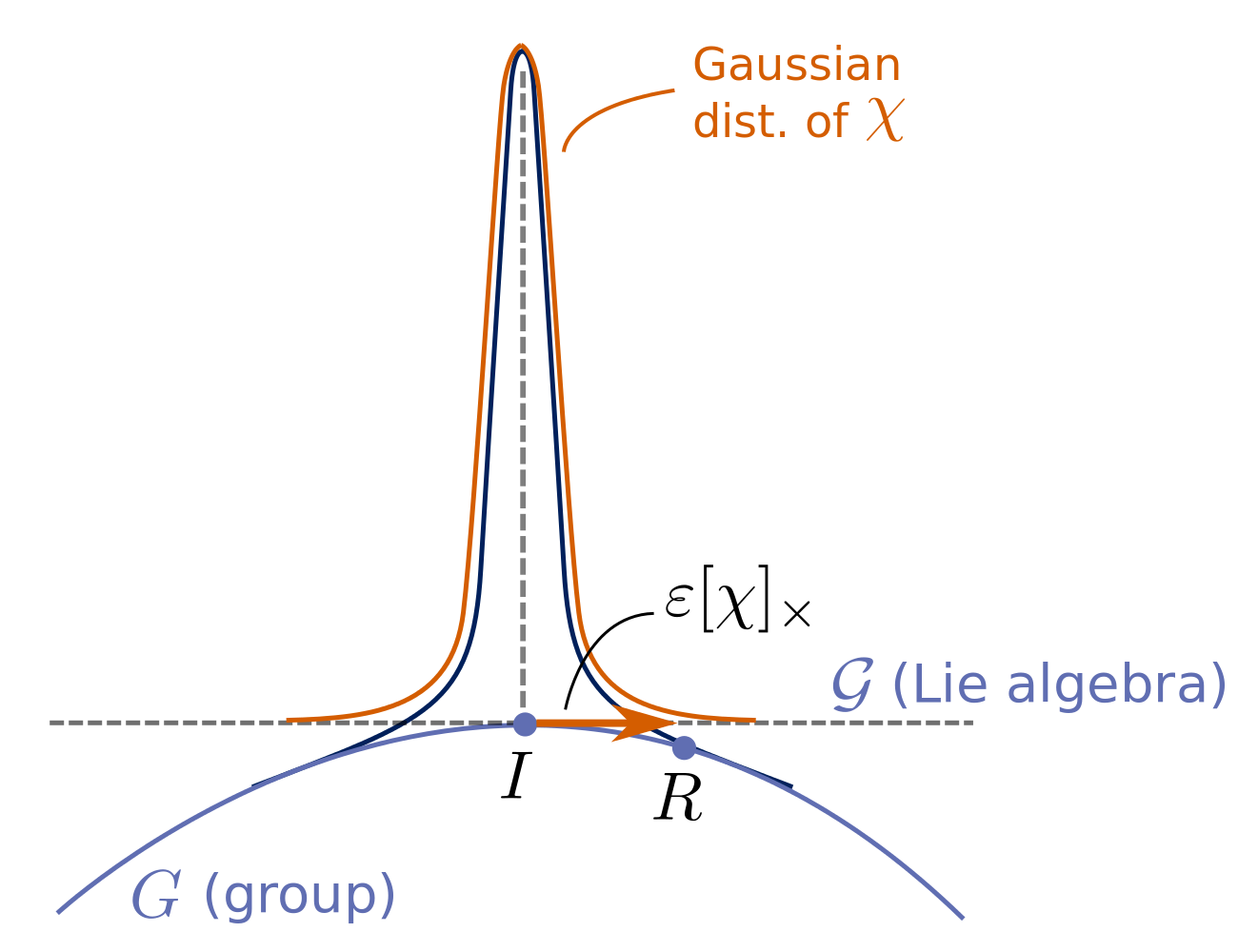}
  \caption{Concentrated distribution on $G$ with mean at the identity, where $R = \exp(\epsilon\,\sk{\dev})$ is close to the mean. 
  The random variable $\dev$ has a Gaussian distribution in $so(3) \cong \Re^3$. 
  A concentrated distribution with non-identity mean is defined using left translation \cite{bourmaud2015EKF}.}
  \label{fig:CGD}
\end{figure}

\medskip

In the attitude estimation problem, the observation model $h(R) = R^T
r$, where $r \in \Re^3$ is a known vector in the inertial
frame (see the accelerometer and magnetometer model
\eqref{eq:accelerometer} and \eqref{eq:magnetometer}).  In this
special case, the constant gain approximation equals the Kalman gain.  

\medskip

\begin{corollary}
 Consider the Poisson equation (33) where the random variable $R = \mu\,\exp(\epsilon\,\sk{\dev})$, 
 where $\dev \in \Re^3$ is a Gaussian random variable with mean $0$ and covariance $\cov$. 
 Let $h_j(R) = e_j^T R^T r$ for $j=1,2,3$. 
 Then, in the asymptotic limit as $\epsilon \rightarrow 0$, the gain function $\K$ is given by,
 \begin{equation}
  \K = \frac{1}{\sigmaW^2} \, \cov H^T + O(\epsilon), 
  \label{eq:constant_gain_SO3}
 \end{equation}
 where $H := \sk{\mu^T r}$. 
 \qed
 \label{cor:constant_gain_SO3}
\end{corollary}

\medskip

\subsection{FPF with Concentrated Distributions}
\label{sec:FPF_EKF}  

Consider the attitude estimation problem, 
\begin{align}
 \dR_t & = R_t\,\Omega_t\,\dt + R_t\circ \sk{\epsilon\,\sigmaB\,\dB_t}, \\
 \dZ_t & = R_t^T r\,\dt + \epsilon\,\sigmaW\,\dW_t, 
\end{align}
with initial condition $R_0 = \mu_0\,\exp(\epsilon\,\sk{\dev_0})$, 
where $\mu_0$ is exactly known and $\dev_0 \sim \clN(0, \cov_0)$. 

The FPF for attitude estimation is given by \eqref{eq:particle_dyn_SO3} 
where the gain is obtained by solving the Poisson equation \eqref{eq:BVP_SO3}. 
For small $\epsilon$ and small time $t \in [0, \,\epsilon T]$, 
the constant gain approximation is used based on Corollary \ref{cor:constant_gain_SO3}. 
The resulting FPF is then given by, 
\begin{equation}
 \dR_t^i = R_t^i\,\Omega_t\,\dt + R_t^i\circ\sk{\epsilon\,\sigmaB\,\dB_t^i} + R_t^i\,\sk{\K_t\circ\dI_t^i}, 
 \label{eq:FPF_concentrated}
\end{equation}
where $\K_t := \frac{1}{\sigmaW^2}\, \cov_t \, H_t^T$ is the constant gain, and $H_t = \sk{\mu_t^T r}$. 

In the following theorem, it is shown that $\mu_t$ and $\cov_t$ evolve according to the 
equations that are closely related to the left invariant EKF. 
The proof is contained in Appendix \ref{apdx:mean_cov_new}. 

\medskip

\begin{theorem}
 Consider the FPF \eqref{eq:FPF_concentrated} where $\K_t$ is given by the constant gain approximation. 
 Suppose that over a time horizon, $ R_t^i = \mu_t\,\exp(\epsilon\,\sk{\dev_t^i})$ where $\dev_t^i \sim \clN(0, \cov_t)$. 
 Then, in the asymptotic limit as $\epsilon \rightarrow 0$,
 $\mu_t$ and $\cov_t$ evolve according to the respective sdes, 
 \begin{align}
  \ud\mu_t = & ~\mu_t \, \Omega_t \, \dt + \mu_t \, \sk{\,\K_t \circ \ud \innov_t\,}, \label{eq:d_mu} \\
  \ud \cov_t = & ~(A_t\dt - \sk{\K_t \dI_t})\,\cov_t + \cov_t\,(A_t\dt - \sk{\K_t \dI_t})^T + \sigmaB^2 \, I\,\dt 
     - \frac{1}{\sigmaW^2}\,\cov_t H_t^T H_t \cov_t\,\dt, \label{eq:d_sigma}
 \end{align}
 where $A_t = -\Omega_t$ and $\ud \innov_t = \ud Z_t - \mu_t^T r \, \dt$.
 \qed
 \label{thm:mean_cov}
\end{theorem}

\medskip 

The equation for the mean \eqref{eq:d_mu} is identical to the left invariant EKF~\cite{bonnabel2009IEKF}. 
The equation of the covariance \eqref{eq:d_sigma} includes additional
terms that depend on the innovation process $\innov_t$.
Analogous stochastic terms for updating the covariance, though in a discrete-time setting, have also appeared in \cite{bourmaud2015EKF}, 
where these terms are induced by the re-parametrization step in the observation update. 
Related results on error propagation and Bayesian fusion in matrix Lie
groups also appear in~\cite{wang2006error, chirikjian2014cdc, wolfe2011bayesian}.

\section{Numerics}
\label{sec:simulations}

In this section, results of two numerical studies are presented for
filters on $SO(3)$: (i) in \Sec{sec:simulation_SO3}, an attitude estimation problem; 
and (ii) in \Sec{sec:simulation_SO2}, a filtering problem for 
a bimodal prior distribution supported on a subgroup $SO(2)$.

\subsection{Attitude Estimation}
\label{sec:simulation_SO3}

Consider an attitude estimation problem with observations from both
accelerometer and magnetometer,
\begin{subequations}
  \begin{align}
    \ud q_t & = \frac{1}{2} \, q_t \otimes \big( \omega_t \ud t + \sigma_B \, \ud B_t \big), \label{eq:dyn_simulation} \\
    \ud Z_t & = 
    \begin{bmatrix}
    -R(q_t)^T & 0 \\
    0 & R(q_t)^T
    \end{bmatrix}
    \begin{bmatrix}
    r^g \\
    r^b
    \end{bmatrix} \ud t + 
    \sigma_W \, \ud W_t, \label{eq:obs_simulation}
  \end{align}
\end{subequations}
where the model for angular velocity is taken from~\cite{zamani2013thesis},
\begin{equation*}
 \omega_t = \Big( \sin \bpar{\frac{2\pi}{15}t},~-\sin \bpar{\frac{2\pi}{18}t+\frac{\pi}{20}},~\cos \bpar{\frac{2\pi}{17}t} \Big),
\end{equation*}
and $r^g=(0,0,1)$, $r^b=(1/\sqrt{2},0,1/\sqrt{2})$ are assumed to be aligned with the gravity 
and the local magnetic field, respectively.

The following attitude filters are simulated for the comparison:
\begin{enumerate}
 \item MEKF:  the multiplicative EKF algorithm described in \cite{markley2003attitude, trawny2005indirect} 
   using the modified Rodrigues parameter. 
 \item USQUE: the unscented quaternion estimator described in
   \cite{crassidis2003unscented} also using the modified Rodrigues parameter.
 \item LIEKF: the left invariant EKF algorithm described in
   \cite{bonnabel2009IEKF}. The equations for the conditional mean and
   covariance are also discussed in Sec. \ref{sec:FPF_EKF}.
 \item IEnKF: the invariant ensemble Kalman filter described in \cite{barrau2015TAC}. 
 \item FPF-G: the FPF using the Galerkin gain function approximation:
   Algorithm~2 in Sec. \ref{sec:galerkin_gain} with the nine basis functions in Table~\ref{table:wigner_d} in Appendix \ref{apdx:basis_SO3}.
 \item FPF-K: the FPF using the kernel-based gain function approximation: Algorithm~3 in Sec. \ref{sec:kernel_gain} with the parameter $\epsilon=1$.
 \item FPF-C: the FPF using the constant gain approximation described in \Sec{sec:CGA}.
\end{enumerate}

The performance metric is evaluated in terms of the {\em rotation
angle error} defined as follows: Let $q_t$ and $\hat{q}_t$ denote
the true and estimated attitude, respectively, at time $t$.  The
estimation error is defined as $\delta q_t := \hat{q}_t^{-1}\otimes
q_t$  and the rotation angle error $\delta\phi_t := 2 \arccos(|\delta
q_t^0|) \in [0,\pi]$, where $\delta q_t^0$ is the first component of $\delta q_t$. 

In an experiment, each filter is simulated over $M$ independent Monte Carlo runs. 
For the $j$-th Monte Carlo run, $\delta\phi_t^j$ denotes the rotation angle error 
as a function of time. The {time-averaged error} for the $j$-th run is defined as, 
\begin{equation}
 \langle \delta\phi^j \rangle_T := \frac{1}{T} \int_0^T \delta\phi_t^j \, \dt, 
 \label{eq:ta_j}
\end{equation}
and the time-averaged error of the $M$ runs is defined as, 
\begin{equation}
 \langle \widehat{\delta \phi} \rangle_T 
     := \frac{1}{M} \sum_{j=1}^M \langle \delta\phi_t^j \rangle_T. 
 \label{eq:ta_error}
\end{equation}
The average error of the $M$ Monte Carlo runs as a function of time is defined according to, 
\begin{equation}
 \widehat{\delta\phi}_t := \frac{1}{M} \sum_{j=1}^M \delta\phi_t^j, 
 \label{eq:average_error}
\end{equation}

The simulation parameters are as follows: The simulations are carried
out over a finite time-horizon $t\in[0,T]$ with fixed time step
$\Delta t$.  The filters are all initialized with a Gaussian distribution, 
denoted as $\clN(q_0, \Sigma_0)$, with mean $q_0$  
and $\Sigma_0 = \sigma_0^2 \, I$ is a diagonal matrix representing the
variance in each axis of the Lie algebra.  For the FPF
implementation, the initial set of particles are sampled from this distribution as follows: 
First, $\{v^i\}_{i=1}^N$ are sampled i.i.d. from the Gaussian
distribution $\clN(0,\Sigma_0)$ in $\Re^3$.  Next, the particles
$\{q_0^i\}_{i=1}^N$ are obtained by,
\begin{equation*}
  q_0^i = q_0 \otimes 
  \begin{bmatrix}
  \cos \big( |\nu^i|/2 \big) \\
  \frac{\nu^i}{|\nu^i|}\, \sin \big( |\nu^i|/2 \big)
  \end{bmatrix}.
\end{equation*}
The IEnKF also uses the same number of particles
as the FPF.

The MEKF, the USQUE and the IEnKF are all discrete-time filters.  They
require a discrete-time filtering model that is chosen to be
consistent with the continuous-time model~\eqref{eq:dyn_simulation}-\eqref{eq:obs_simulation}.  For the
discrete-time filters, the sampled observations, denoted as
$\{\tn{Y}\}$, are made at discrete times $\{t_n\}$, whose model is
formally expressed as $\tn{Y} := \frac{\Delta \tn{Z}}{\Dt} = h(\tn{q})
+ \tn{W}^\Delta$ where $\{\tn{W}^\Delta\}$ are i.i.d. with the
distribution $\clN(0,\frac{\sigma_W^2}{\Dt}\,I)$.  Such a model leads
to the correct scaling between the continuous and the discrete-time
filter implementations.

In numerical simulations, it was observed that the continuous-time
filters, especially the FPF-G, are susceptible to numerical
instabilities due to high gain during the initial transients.  The
instability in FPF-G is exacerbated by possible ill-conditioning of
the matrix $A$ in constructing the Galerkin approximation (see Algorithm \ref{alg:galerkin}). 
In order to mitigate the numerical issues observed during the implementation of
the FPF-G algorithm, the discrete time-step during the initial
transients is further sub-divided.  Specifically, for $t < T_f$, the
time interval $[\,t, t+\Dt\,]$ is uniformly divided into $N_f$ 
sub-intervals.  The update step in the FPF (specifically Line 3 --
Line 11 in Algorithm 4) is implemented on each sub-interval by
replacing $\Delta Z_t$ with $\frac{\Delta Z_t}{N_f}$ and $\Dt$ with
$\frac{\Dt}{N_f}$. To provide a fair comparison, the same set of
observations are used for all the continuous-time and the
discrete-time algorithms.   

The nominal parameter values are chosen as: 
$T=3$, $\Dt=0.01$, $N=100$, $M=100$, $T_f=0.2$, $N_f=100$. 
The choice of $T_f$ and $N_f$ may vary according to the severity of numerical issues encountered in practice.  

\begin{figure*}
  \vspace*{-15pt}
  \centering
  \begin{subfigure}{0.44\textwidth}
   \captionsetup{skip=-2pt}
   \includegraphics[width=\textwidth]{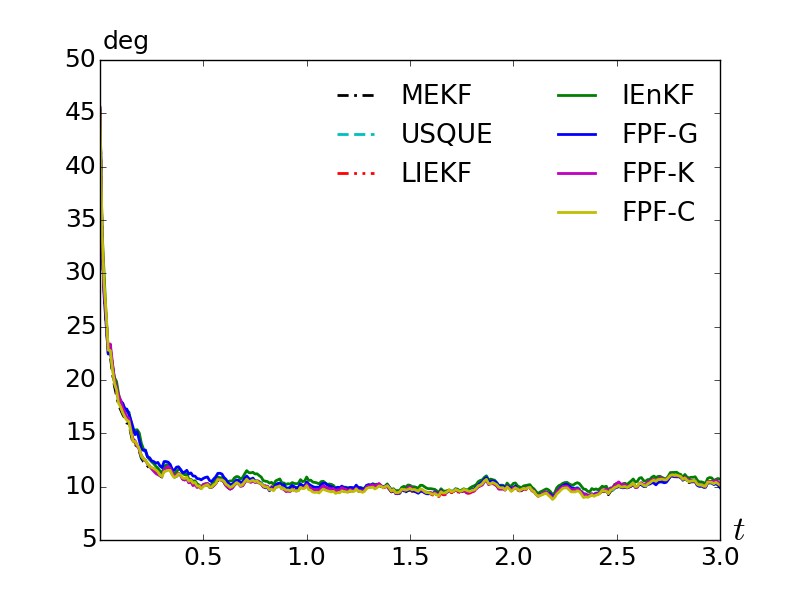}
   \caption{Initial distribution with $\sigma_0=\degg{30}$}
  \end{subfigure}\hspace*{1em}
  \begin{subfigure}{0.44\textwidth}
   \captionsetup{skip=-2pt}
   \includegraphics[width=\textwidth]{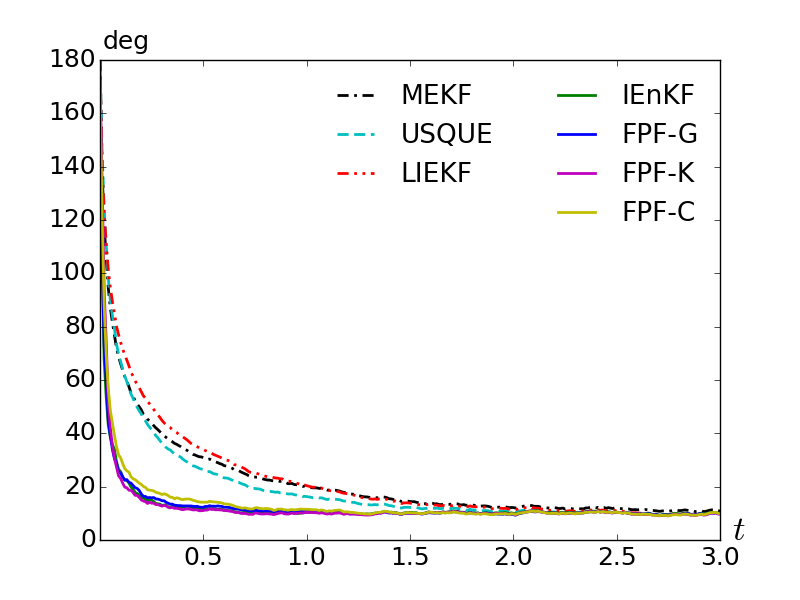}
   \caption{Initial distribution with $\sigma_0=\degg{60}$}
  \end{subfigure}
  \vspace*{-1pt}
  \caption{Comparison of the average error $\widehat{\delta\phi}_t$.}
  \label{fig:init_var}
\end{figure*}

\begin{figure*}
  \vspace*{-10pt}
  \centering
  \begin{subfigure}{0.44\textwidth}
   \captionsetup{skip=-2pt}
   \includegraphics[width=\textwidth]{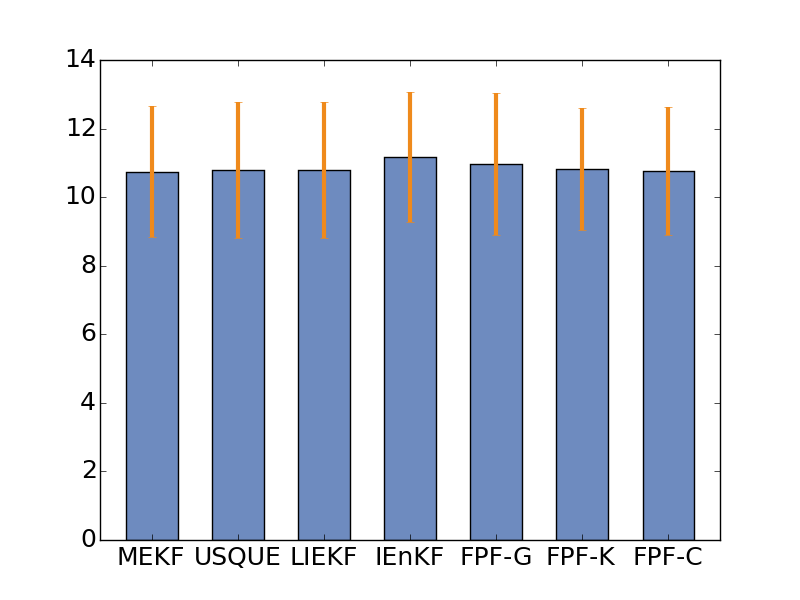}
   \caption{Initial distribution with $\sigma_0=\degg{30}$}
  \end{subfigure}\hspace*{1em}
  \begin{subfigure}{0.44\textwidth}
   \captionsetup{skip=-2pt}
   \includegraphics[width=\textwidth]{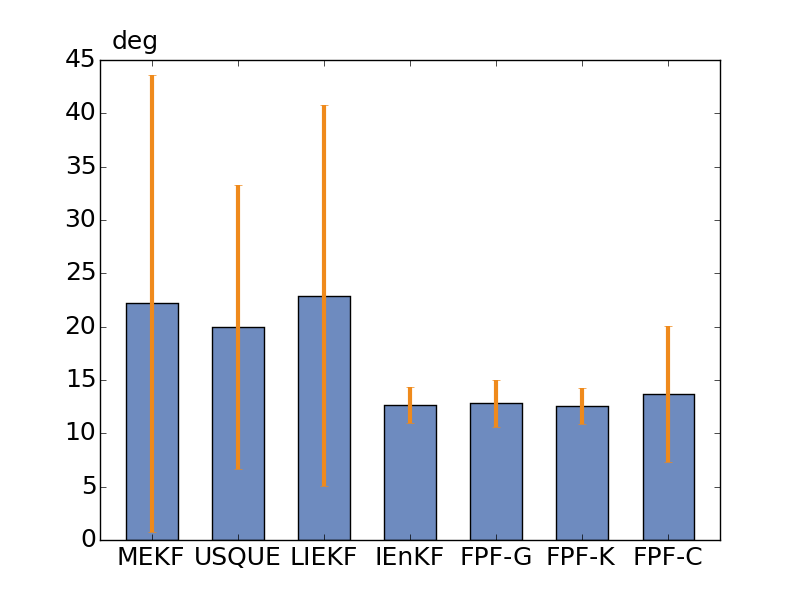}
   \caption{Initial distribution with $\sigma_0=\degg{60}$}
  \end{subfigure}
  \vspace*{-1pt}
  \caption{Statistical analysis of filter performance: The bars indicate the mean and the lines 
  indicate the $\pm 1$ standard deviation of $\{\langle \delta\phi^j \rangle_T\}_{j=1}^M$
  across $M=100$ Monte-Carlo runs}
  \label{fig:variability}
\end{figure*}

\medskip

\begin{figure*}
  \vspace*{-10pt}
  \centering
  \begin{subfigure}{0.44\textwidth}
   \captionsetup{skip=-2pt}
   \includegraphics[width=\textwidth]{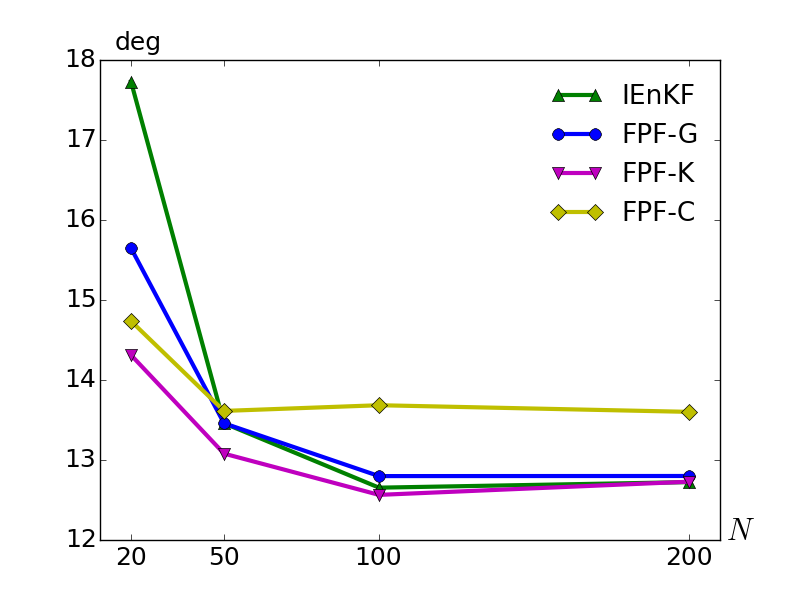}
   \caption{}
  \end{subfigure}\hspace*{1em}
  \begin{subfigure}{0.44\textwidth}
   \captionsetup{skip=-2pt}
   \includegraphics[width=\textwidth]{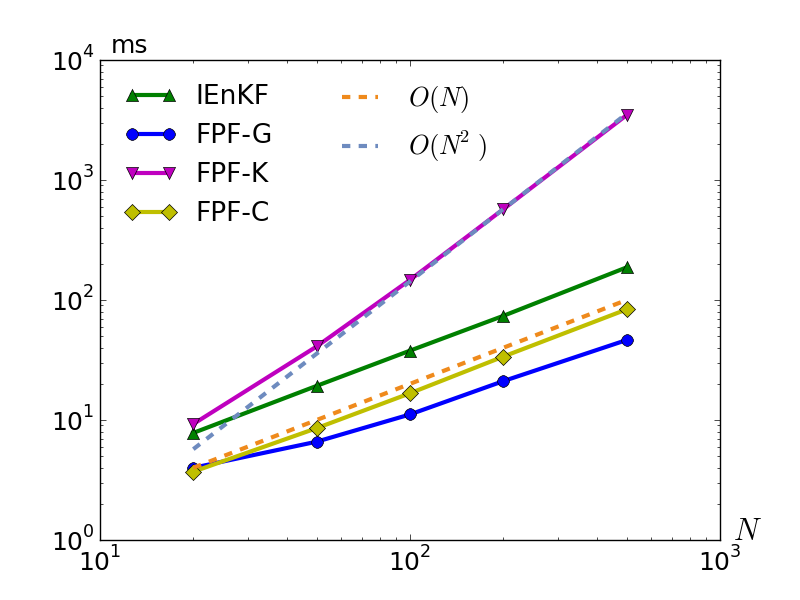}
   \caption{}
  \end{subfigure}
  \vspace*{-1pt}
  \caption{(a): Time-averaged error $\langle \widehat{\delta \phi} \rangle_T$, 
  and (b): mean computational time of a single propagation-update step, both as a function of the number of particles $N$.}
  \label{fig:N_effect}
\end{figure*}

The simulation results are discussed next:

\begin{enumerate}
\item{\bf The average error $\widehat{\delta\phi}_t$ as a function of the initial uncertainty:} 
  Figure~\ref{fig:init_var} depicts the average error $\widehat{\delta\phi}_t$ 
  (see \eqref{eq:average_error}) 
  of the filters over $M=100$ simulation runs, with 
  two choices of initial variance: (a) $\Sigma_0 = 0.5236^2\id$ and
  (b) $\Sigma_0 = 1.0472^2\id$.  The two cases correspond to a
  standard deviation of  $\degg{30}$ and $\degg{60}$, respectively.
  For the two priors, the mean is the same, given by identity quaternion $q_I=(1,0,0,0)$. 
  For case~(a), the target is initialized by sampling from the prior
  distribution.  For case~(b), the target is initialized with a fixed
  attitude\,--\,rotation of $\degg{180}$ about the axis $(3,1,4)$.
  These parameters indicate large estimation error initially for
  case~(b).  

  The results depicted in Figure~\ref{fig:init_var} show that the
  performance is nearly identical across filters for case~(a) when the initial
  uncertainty is small.  For case~(b) when the initial uncertainty is
  large, the particle-based filters exhibit superior performance compared to
  the Kalman filters and unscented filter.  The differences are exhibited
  in the speed of convergence of the estimate to the target with
  the particle filters converging quickly compared to the Kalman filters
  and the unscented filter.  

  As the results in Figure~\ref{fig:init_var} are averaged over multiple
  Monte-Carlo runs, statistical analysis was also carried out to assess
  the variability in performance across runs.  The results of this
  analysis are presented in Figure~\ref{fig:variability}, 
  which depicts the mean and standard deviation of $\{\langle \delta\phi^j \rangle_T\}_{j=1}^M$ 
  (see \eqref{eq:ta_j}). 
  Apart from poorer performance
  on average, the Kalman filters also exhibit a greater variability in
  performance across the Monte-Carlo runs.  
  For some trajectories, the Kalman filters exhibit slow convergence because the gain becomes very small.

\begin{figure*}
  \vspace*{-15pt}
  \centering
  \begin{subfigure}{0.44\textwidth}
   \captionsetup{skip=-2pt}
   \includegraphics[width=\textwidth]{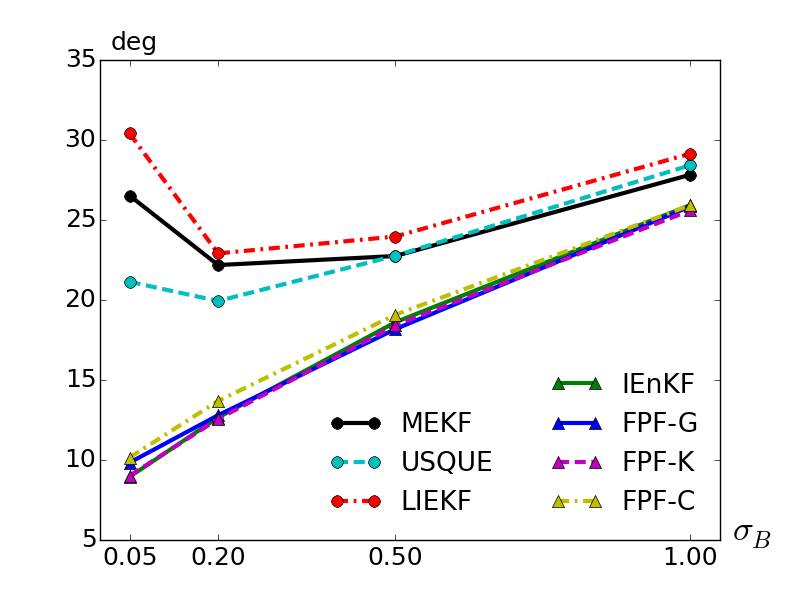}
   \caption{}
  \end{subfigure}\hspace*{1em}
  \begin{subfigure}{0.44\textwidth}
   \captionsetup{skip=-2pt}
   \includegraphics[width=\textwidth]{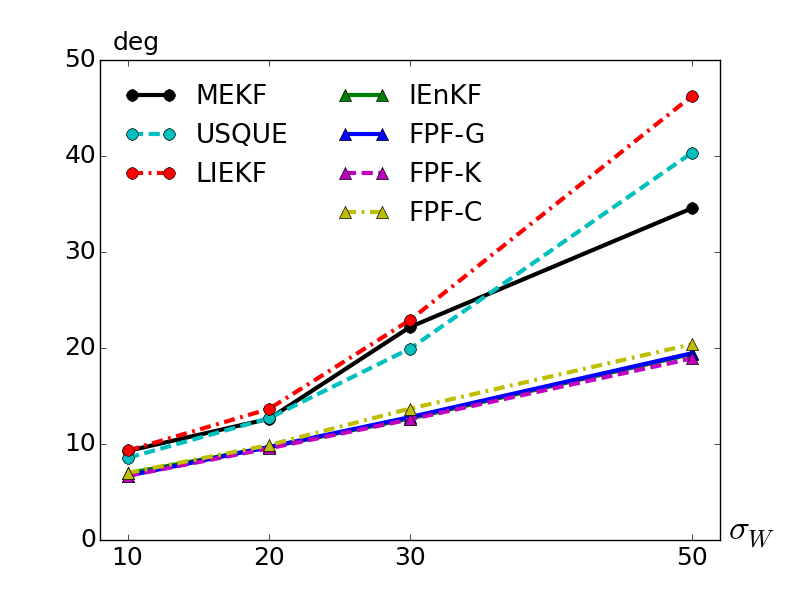}
   \caption{}
  \end{subfigure}
  \vspace*{-1pt}
  \caption{Time-averaged error $\langle \widehat{\delta \phi} \rangle_T$ of filters as a 
  function of $\sigma_B$ and $\sigma_W$. The value of $\sigma_W$ is converted to the standard deviation (in degree) 
  of the corresponding discrete-time observation model.}
  \label{fig:sigma_B_W}
\end{figure*}

\item{\bf The time-averaged error $\langle \widehat{\delta \phi} \rangle_T$ 
  as a function of the process noise:} In this
  simulation, the process noise $\sigma_B \in \{0.05,\, 0.2,\,
  0.5,\,1.0\}$ for fixed $\sigma_W=0.05236$ and prior distribution according to
  case~(b) in Figure \ref{fig:init_var}.  Figure~\ref{fig:sigma_B_W}~(a) depicts
  the time-averaged error $\langle \widehat{\delta \phi} \rangle_T$ 
  (see \eqref{eq:ta_error}) across filters 
  as the process noise parameter is varied. 
  One would have expected the error to increase
  monontonically with the $\sigma_B$ value.  The fact that such is not
  the case for the Kalman filters indicates that the relatively poor
  performance of the Kalman filters for small values of process noise
  is an artifact of the linearization assumption that leads to overly small
  gains. These small gains adversely effect the filter performance 
  during the initial transients. 

\item{\bf The time-averaged error $\langle \widehat{\delta \phi} \rangle_T$ 
  as a function of the observation noise:} In this
  simulation, the observation noise parameter $\sigma_W \in \{0.01745,\,0.03491,\,0.05236,\,0.08727\}$ 
  for fixed $\sigma_B=0.2$ and prior distribution according to
  case~(b) in Figure \ref{fig:init_var}. The $\sigma_W$ parameter values
  correspond to the choice of the standard deviation of $\degg{10},
  \,\degg{20}, \,\degg{30}$ and $\degg{50}$ in the discrete-time
  model. Figure~\ref{fig:sigma_B_W}~(b) depicts the time-averaged error $\langle \widehat{\delta \phi} \rangle_T$. 
  As expected, the error deteriorates as the observation noise increases.
  The particle filters not only continue to exhibit better performance
  but also the performance deterioration is more graceful for larger values of $\sigma_W$.

\item{\bf The time-averaged error $\langle \widehat{\delta \phi} \rangle_T$ as a function of $N$:} 
In this simulation, $N \in \{20, 50, 100, 200\}$ in the particle filters,  
for a fixed $\sigma_B=0.2$, $\sigma_W=0.05236$, and prior distribution according to case (b) in Figure \ref{fig:init_var}. 
Figure~\ref{fig:N_effect}~(a) depicts the time-averaged error $\langle \widehat{\delta \phi} \rangle_T$. 
For all the particle filters, $N=50$ particles is seen to be sufficient. 
For fewer than $50$ particles, 
the FPF-G and the IEnKF exhibit performance deterioration  
as insufficient number of particles leads to issues in the gain computation. 
Numerically, the FPF-K is seen to be the best algorithm for small value of $N$. 

\item{\bf Computational times as a function of $N$:} In this simulation,
$N \in \{20, 50, 100, 200, 500\}$. The mean computational
time (per propagation-update step of the algorithm, averaged over 100 Monte Carlo runs) is depicted
as a function of $N$ in Figure~\ref{fig:N_effect} (b).  The
$O(N)$ and $O(N^2)$ lines are included to aid the comparison. 
The computational cost of particle filters scale linearly with $N$
except the kernel method which scales quadratically. 
For online computations, both FPF-G and FPF-C have lower computational burden compared to IEnKF. 
However, for the IEnKF algorithm, the gain computation \,--\,
which contributes to most of the computation load \,--\, can be
implemented offline \cite{barrau2015TAC}.  
The experiments were conducted on a platform with an Intel i3-2120 3.3GHz CPU. 

\end{enumerate}

\subsection{Filtering with a Bimodal Distribution}
\label{sec:simulation_SO2}

\begin{figure*}
  \centering
    \includegraphics[width=0.9\textwidth]{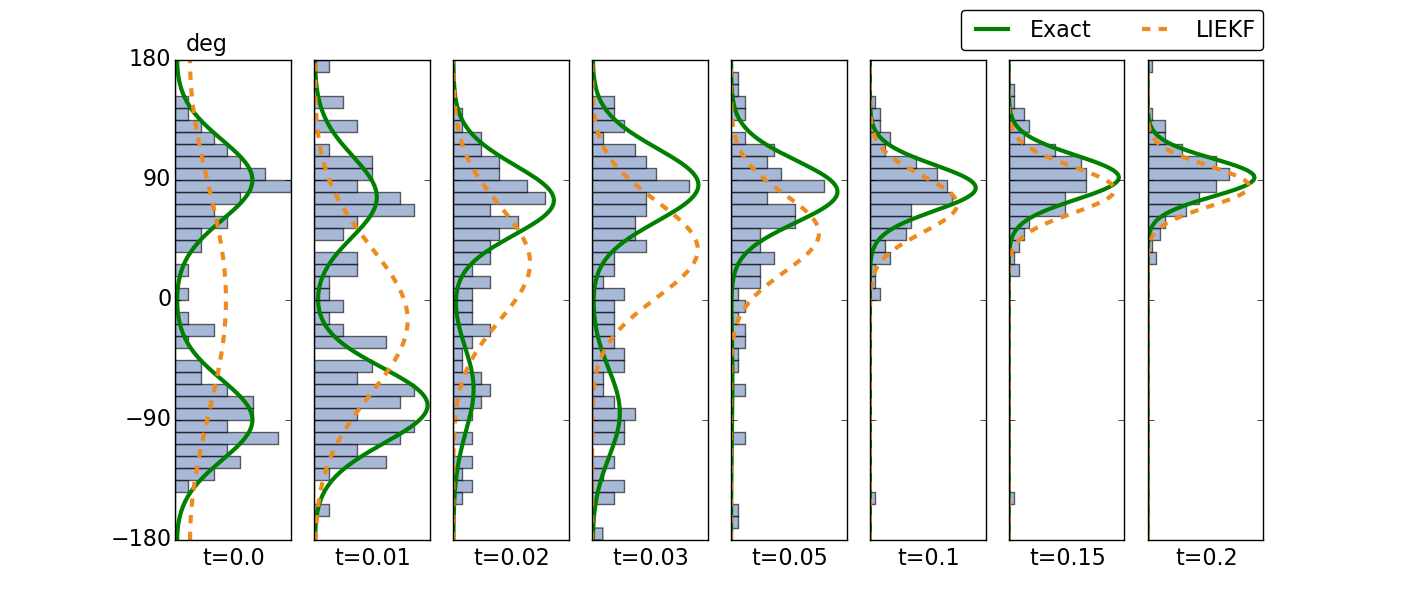}
  \vspace*{-5pt}
  \caption{Density evolution on $SO(2)$ with bimodal distribution.}
  \label{fig:density_evolution}
\end{figure*}

In this section, we consider the following static model in $SO(3)$: 
  $$ \dq_t = \frac{1}{2}\,q_t \otimes \omega_t\,\dt, $$
where $\omega_t = (0,0,0)$. The prior distribution 
is assumed to be supported on the subgroup $SO(2)$, parametrized by the angle $\theta \in [-\pi,\pi)$. 
Its density is denoted as $\rho^*_0(\theta)$. 
An arbitrary element in $SO(2)$ is represented as 
$q = \bpar{\cos(\frac{\theta}{2}),\,0,\,0,\,\sin(\frac{\theta}{2})}$. 

The observation model is of the following form:
  $$ \dZ_t = h(\theta_t)\,\dt + \sigma_W\,\dW_t, $$
where $h(\theta) = \bpar{\cos(\theta), -\sin(\theta)}$, and $W_t$ is a standard Wiener process in $\Re^2$. 

Since the process is static, the density of the posterior distribution has a closed-form Bayes' formula: 
\begin{equation}
 \rho^*(\theta, t) = \text{(const.)}~\exp\Bpar{\frac{1}{\sigma_W^2}\,h^T(\theta)\,Z_t 
 - \frac{1}{2\sigma_W^2}\,|h(\theta)|^2\,t} \, \rho^*_0(\theta). 
 \label{eq:exact_dist}
\end{equation}

For the numerical results described next, the FPF is simulated according to \eqref{eq:particle_dyn_SO3}:
  $$ \dq_t^i = \frac{1}{2}\,q_t^i \otimes \Big[ \K(q_t^i,t)\circ\bpar{\dZ_t - \frac{h(q_t^i) + \hat{h}}{2}\,\dt} \Big], $$
where $q_0^i$ are sampled i.i.d. from the prior $\rho^*_0$. 

The simulation parameters are as follows: The prior is a mixture of two Gaussians, 
$\clN\bpar{-\mu_0, \sigma_0^2}$ and $\clN\bpar{\mu_0, \sigma_0^2}$, 
with equal weights, where $\mu_0 = \degg{90}$ and $\sigma_0=\degg{30}$. 
The observation noise parameter $\sigma_W=0.12$, 
and the unknown state is initialized as $q_0 = (1/\sqrt{2},\,0,\,0,\,1/\sqrt{2})$, 
which corresponds to $\theta_0 = \degg{90}$. 
The simulations are carried out over $t \in [0, 0.2]$ with fixed time step $\Dt=0.01$. 
The FPF-K with $N=100$ and $\epsilon=0.2$ is simulated, together with the The LIEKF for a comparison. 

Figure \ref{fig:density_evolution} depicts the simulation results including 
the exact posterior (see \eqref{eq:exact_dist}), the histogram of the particles, and the LIEKF solution. 
This example shows that the FPF algorithm can easily handle a general class of non-Gaussian distributions.

\appendices

\section{Proof of Proposition \ref{prop:evolution}}
\label{apdx:prop_evolution}

For any function $f \in \Cinfc(G)$, $f(X_t^i)$ is a continuous semimartingale that satisfies \cite{hsu2002}, 
\begin{align}
 \ud f(X_t^i) = & ~(V_0+\u)\cdot f(X_t^i)\ud t + V_\ip\cdot f(X_t^i) \circ \ud B_t^{\ip,i} \notag \\ 
 & ~+ \K_\io\cdot f(X_t^i) \circ \ud Z_t^\io.
 \label{eq:df_strat}
\end{align}

For the ease of taking the expectation,
we convert \eqref{eq:df_strat} to its It\^{o} form (see Theorem 1.2 in \cite{watanabe1981}):
For real-valued continuous semi-martingales $A, B, C$,
\begin{align}
 & A\circ \ud B = A \ud B + \frac{1}{2}\ud A \ud B, \label{eq:eq1} \\
 & (A\circ\ud B) \ud C = A (\ud B \ud C). \label{eq:eq2}
\end{align}

For the second term on the right hand side of \eqref{eq:df_strat}, taking $A$ in \eqref{eq:eq1} to be 
$V_\ip\cdot f(X_t^i)$ and $B$ to be $B_t^{\ip,i}$,
\begin{equation}
 V_\ip\cdot f(X_t^i) \circ \ud B_t^{\ip,i} = V_\ip\cdot f(X_t^i) \ud B_t^{\ip,i}
 + \frac{1}{2} \ud (V_\ip\cdot f)(X_t^i) \ud B_t^{\ip,i}.
 \label{eq:convert_E}
\end{equation}
Replacing $f$ by $V_\ip\cdot f$ in \eqref{eq:df_strat},
\begin{align*}
 \ud (V_\ip\cdot f) = & ~(V_0+\u)\cdot (V_\ip\cdot f)\ud t + 
 V_\ipp\cdot(V_\ip\cdot f) \circ \ud B_t^{\ipp,i} + \K_\io\cdot (V_\ip\cdot f) \circ \ud Z_t^\io.
\end{align*}
Using \eqref{eq:eq2} and It\^{o}'s rule 
($\dB_t^{\ip,i}\dB_t^{\ipp,i} = \delta_{\ip,\ipp}\,\dt$, 
$\dB_t^{\ip,i}\dt = 0$, and $\dB_t^{\ip,i}\dZ_t^\io = 0$ for all $\ip, \beta, \io$), 
\begin{equation*}
 \ud (V_\ip\cdot f)(X_t^i) \ud B_t^{\ip,i} = \sum_{\ip=1}^{r} V_\ip\cdot(V_\ip\cdot f)(X_t^i) \ud t,
\end{equation*}
which when substituted in \eqref{eq:convert_E} yields,
 $$ V_\ip\cdot f(X_t^i) \circ \ud B_t^{\ip,i} = V_\ip \cdot f(X_t^i) \ud B_t^{\ip,i} 
  + \frac{1}{2} \sum_{\ip=1}^r V_\ip\cdot(V_\ip\cdot f)(X_t^i) \ud t. $$

The third term on the right hand side of \eqref{eq:df_strat} is similarly converted. 
The It\^{o} form of \eqref{eq:df_strat} is then given by,
\begin{equation*}
 \ud f(X_t^i) = \mathcal{L}f(X_t^i) \, \ud t + V_\ip \cdot f(X_t^i) \ud B_t^{\ip,i}
 + \K_\io\cdot f(X_t^i,t) \ud Z_t^\io,
\end{equation*}
where the operator $\mathcal{L}$ is defined by,
\begin{equation*}
 \mathcal{L}f := (V_0+\u)\cdot f + \frac{1}{2} \sum_{\ip=1}^r V_\ip \cdot (V_\ip\cdot f) 
 + \frac{1}{2} \sum_{\io=1}^m \K_\io \cdot (\K_\io\cdot f).
\end{equation*}

In its integral form, 
\begin{align}
 f(X_t^i) = & ~f(X_0^i) + \int_0^t \mathcal{L}f(X_s^i)\ud s 
            + \int_0^t V_\ip \cdot f(X_s^i) \ud B_s^{\ip,i} + \int_0^t \K_\io \cdot f(X_s^i) \ud Z_s^\io. \notag
\end{align}
By taking conditional expectation on both sides, interchanging expectation and integration 
(see Lemma 5.4 in \cite{xiong2008introduction})
and noting the fact that $B_t^{\ip,i}$ is a Wiener process, 

 $$ \pi_t (f) = \pi_0(f) 
   + \int_0^t \pi_s( \mathcal{L}f ) \ud s 
   + \int_0^t \pi_s( \K_\io \cdot f ) \ud Z_s^\io,
 $$
which is the desired formula \eqref{eq:forward}.

\begin{table*}[t] 
 \caption{{Basis functions on $SO(3)$}}
 \centering
 \begin{tabular}{c||c|c||c|c|c}
  \hline
             & expression in $R$ & expression in $q$  & $E_1\cdot$    & $E_2\cdot$            & $E_3\cdot$ \\ \hline
    $\psi_1$ & $R_{33}$  & $2(q_0^2+q_3^2)-1$  & $2(-q_0q_1-q_2q_3)$  & $2(-q_0q_2+q_1q_3)$   & $0$        \\ \hline
    $\psi_2$ & $R_{13}$  & $2(q_0q_2+q_1q_3)$  & $2(q_0q_3-q_1q_2)$   & $2(q_0^2+q_1^2)-1$    & $0$        \\ \hline
    $\psi_3$ & $-R_{23}$ & $2(q_0q_1-q_2q_3)$  & $2(q_0^2+q_2^2)-1$   & $2(-q_0q_3-q_1q_2)$   & $0$        \\ \hline
    $\psi_4$ & $R_{31}$  & $2(-q_0q_2+q_1q_3)$ & $0$                  & $-2(q_0^2+q_3^2)+1$   & $2(q_0q_1+q_2q_3)$   \\ \hline
    $\psi_5$ & $R_{32}$  & $2(q_0q_1+q_2q_3)$  & $2(q_0^2+q_3^2)-1$   & $0$                   & $2(q_0q_2-q_1q_3)$   \\ \hline
    $\psi_6$ & $(1/2)(R_{21}-R_{12})$ & $2q_0q_3$     & $-q_0q_2-q_1q_3$     & $q_0q_1-q_2q_3$       & $q_0^2-q_3^2$ \\ \hline
    $\psi_7$ & $(1/2)(R_{11}+R_{22})$ & $q_0^2-q_3^2$ & $-q_0q_1+q_2q_3$     & $-q_0q_2-q_1q_3$      & $-2q_0q_3$    \\ \hline
    $\psi_8$ & $(1/2)(R_{21}+R_{12})$ & $2q_1q_2$     & $q_0q_2+q_1q_3$      & $q_0q_1-q_2q_3$       & $q_2^2-q_1^2$ \\ \hline
    $\psi_9$ & $(1/2)(R_{11}-R_{22})$ & $q_1^2-q_2^2$ & $q_0q_1-q_2q_3$      & $-q_0q_2-q_1q_3$      & $2q_1q_2$     \\
  \hline
 \end{tabular}
 \label{table:wigner_d}
\end{table*}

\section{Proof of Theorem \ref{thm:consistency}}
\label{apdx:thm_consistency}

Using~\eqref{eq:K-S} and~\eqref{eq:forward} and the expressions for the
operators $\clL^*$ and $\clL$, it suffices to show that,
\begin{align}
 & \pi_s(\u\cdot f) \ud s + \frac{1}{2} \sum_{\io=1}^m \pi_s \big( \K_\io \cdot (\K_\io \cdot f) \big)\,\ds + \pi_s(\K_\io \cdot f) \ud Z_s^\io \notag \\
 & = \sum_{\io=1}^m \big(\pi_s(fh_\io)-\pi_s(h_\io)\pi_s(f)\big) \big(\ud Z_s^\io - \pi_s(h_\io)\ud s \big)
 \label{eq:WTS}
\end{align}
for all $0 \leq s \leq t$ and all $f\in\Cinfc(G)$.

On taking $\psi=f$ in \eqref{eq:BVP} and using the
formula~\eqref{eq:grad_inner_product} for the inner product, 
\begin{equation}
 \pi_s(\K_\io \cdot f) = \pi_s \big( (h_\io-\pi_s(h_\io))f \big).
 \label{eq:Kf}
\end{equation}

Using the expression~\eqref{eq:optimal_u} for the control function and noting that $\hat{h}_\io=\pi_s(h_\io)$,
  $$ \u \cdot f = -\frac{1}{2} \sum_{\io=1}^m \bpar{h_\io-\pi_s(h_\io)} \, \bpar{\K_\io \cdot f} 
  - \sum_{\io=1}^m \pi_s(h_\io) \, \bpar{\K_\io \cdot f}. $$
Using \eqref{eq:Kf} repeatedly then leads to, 
\begin{align}
 \pi_s(\u\cdot f) = & ~-\frac{1}{2} \sum_{\io=1}^m \pi_s \big( \K_\io \cdot(\K_\io \cdot f) \big) - \sum_{\io=1}^m  \pi_s(h_\io) \pi_s \big( (h_\io-\pi_s(h_\io))f \big).
 \label{eq:uf}
\end{align}

The desired equality~\eqref{eq:WTS} is now verified by substituting in \eqref{eq:Kf} and \eqref{eq:uf}.

\section{Basis Functions on $SO(3)$}
\label{apdx:basis_SO3} 

The eigenfunctions of the Laplacian on $SO(3)$ are determined by the matrix elements of the 
irreducible unitary representations of $SO(3)$ (see Sec. 9.4 in \cite{chirikjian2016harmonic}). 
The eigenfunctions associated with the smallest non-zero eigenvalue are tabulated in Table \ref{table:wigner_d}, 
expressed using both the rotation matrix and the quaternion. 
In order to compute the matrix $A$ in the Galerkin gain function approximation, 
the formulae for $E_1\cdot\psi_l,~E_2\cdot\psi_l~\text{and}~E_3\cdot\psi_l$ are also provided, 
where $\{E_1,\,E_2,\,E_3\}$ denote the basis of $so(3)$ given by \eqref{eq:basis_so3}.

\section{Proof of Proposition \ref{prop:constant_gain}}
\label{apdx:constant_gain} 

Using the basis $\{E_1,\,E_2,\,E_3\}$ of the Lie algebra $so(3)$, 
the strong form of the Poisson equation (see \eqref{eq:PE_strong_form}) is expressed as, 
 \begin{equation}
  \sum_{n=1}^3 E_n\cdot\bpar{\rho\,\k_{n,\io}} = -(h_\io - \hat{h}_\io)\,\rho,
  \label{eq:poisson_coord}
 \end{equation}
 where $\rho$ represents the density function associated with the distribution $\pi$, 
 and $\k_{n,\io}$ denotes the $(n,j)$-th element of $\K$. Accordingly, the $\io$-th column of $\K$ 
 is $[\K]_\io = (\k_{1,\io},\,\k_{2,\io},\,\k_{3,\io})$. 
 
 Since $\dev$ is Gaussian, $\rho$ is of the form, 
 \begin{equation*}
  \rho(R) = C(R)\,\exp\Bpar{-\frac{1}{2} \big[\,\log(\mu^T R)\,\big]^{\vee\,T} \Sigma^{-1} \big[\,\log(\mu^T R)\,\big]^{\vee}}, 
 \end{equation*}
 where $C(R) \approx 1/\sqrt{(2\pi)^3|\Sigma|}$ 
 (i.e., a constant) when the distribution is concentrated \cite{wang2006error}, and $\Sigma := \epsilon^2\,\cov$. 
 By direct calculation, we have, 
 \begin{equation}
  E_n\cdot\rho = -\frac{1}{\epsilon}\,\dev^T \inv{\cov} M \, e_n\,\rho, 
  \label{eq:d_rho}
 \end{equation}
where $M = I + O(\epsilon)$. 
 
Using \eqref{eq:d_rho}, the left-hand-side of \eqref{eq:poisson_coord} is expanded as follows, 
\begin{align*}
 \sum_{n=1}^3 E_n\cdot\bpar{\rho\,\k_{n,\io}} & = \sum_{n=1}^3 \bpar{E_n\cdot\rho}\,\k_{n,\io} + \rho\,\sum_{n=1}^3 E_n\cdot\k_{n,\io} 
 = -\frac{1}{\epsilon}\,\dev^T \inv{\cov} M \, [\K]_\io\,\rho + \rho\,\sum_{n=1}^3 E_n\cdot\k_{n,\io}.
\end{align*}

The Taylor expansion of $h_\io$ is given by, 
  $$ h_\io(R) = h_\io(\mu) + \epsilon \, \dev^T H_\io + O(\epsilon^2), $$
where $H_\io := \bpar{ E_1\cdot h_\io(\mu),\,E_2\cdot h_\io(\mu),\,E_3\cdot h_\io(\mu) }$. 
Using the fact that $\pi(\dev) = 0$, we have $\hat{h}_\io = \pi(h_\io) =  h_\io(\mu) + O(\epsilon^2)$, 
leading to $ h_\io - \hat{h}_\io = \epsilon\,\dev^T H_\io + O(\epsilon^2) $. 
Hence, \eqref{eq:poisson_coord} becomes, 
  $$ -\dev^T \inv{\cov} M \, [\K]_\io + \epsilon\,\sum_{n=1}^3E_n\cdot\k_{n,\io}
     = -\frac{1}{\overline{\sigma}_W^2}\,\dev^T H_\io + O(\epsilon). $$ 
In the asymptotic limit as $\epsilon\rightarrow 0$, $[\K]_\io = \frac{1}{\overline{\sigma}_W^2} \, \cov \, H_\io + O(\epsilon)$.

\section{Proof of Theorem \ref{thm:mean_cov}}
\label{apdx:mean_cov_new}

Under the constant gain approximation, the FPF is given by (see \eqref{eq:FPF_concentrated}),
\begin{equation}
 \dR_t^i = R_t^i\,\Omega_t\,\dt + R_t^i\circ\sk{\epsilon\,\sigmaB\,\dB_t^i} + R_t^i\,\sk{\K_t\circ\dI_t^i}, 
 \label{eq:FPF_concentrated_apdx}
\end{equation}
where $\K_t := \frac{1}{\sigmaW^2}\, \cov_t \, H_t^T$ and $H_t = \sk{\mu_t^T r}$. 
A concentrated distribution is assumed, i.e., , 
\begin{equation}
 R_t^i = \mu_t\,\exp(\epsilon\,\sk{\dev_t^i}) = \mu_t + \epsilon\,\mu_t\,\sk{\dev_t^i} + O(\epsilon^2), 
 \label{eq:representation}
\end{equation}
where $\dev_t^i \sim \clN(0, \cov_t)$. 

The evolution equation for the mean $\mu_t$ and covariance $\cov_t$ are derived using a perturbation analysis approach. 
We begin by simplifying the modified form of the innovation error,
\begin{align}
 \dI_t^i & = \dZ_t - \frac{R_t^{i\,T}r + \pi_t(R_t^{i\,T}r)}{2}\dt = \dZ_t - \frac{R_t^{i\,T}r + \mu_t^T r}{2}\dt + O(\epsilon^2) \notag \\ 
         & = \dI_t - \frac{R_t^{i\,T} r - \mu_t^T r}{2}\,\dt + O(\epsilon^2) 
           = \dI_t - \frac{1}{2}\,\epsilon\,H_t\,\dev_t^i\,\dt + O(\epsilon^2), \label{eq:dI_approx}
\end{align}
where $\dI_t = \dZ_t - \mu_t^T r\,\dt$, and we have used the fact that $\pi_t(R_t^{i\,T}r) = \mu_t^T r + O(\epsilon^2)$.  

On substituting \eqref{eq:representation} and \eqref{eq:dI_approx} into the FPF \eqref{eq:FPF_concentrated_apdx} 
and matching terms, the $O(1)$ balance gives, 
\begin{equation}
 \ud \mu_t = \mu_t\,\Omega_t\,\dt + \mu_t \sk{\K_t \circ \dI_t}. 
 \label{eq:d_mu_apdx}
\end{equation}
The $O(\epsilon)$ balance gives, 
\begin{align*}
 \ud (\mu_t \, \sk{\dev_t^i}) = & ~\mu_t\,\sk{\dev_t^i}\,\Omega_t\,\dt + \sigmaB\,\mu_t\,\sk{\dB_t^i} 
 + \mu_t\,\sk{\dev_t^i} \sk{\K_t\circ\dI_t} - \frac{1}{2}\,\mu_t\,\sk{\K_t \, H_t \, \dev_t^i} \dt.
\end{align*}
Using the formula \eqref{eq:d_mu_apdx}, this is simplified to obtain the following equation of $\dev_t^i$, expressed in its \ito~form:  
\begin{align}
 \ud \dev_t^i = & A_t \dev_t^i\, \dt + \sigma_B\,\dB_t^i 
    - \frac{1}{2}\,\K_t\,H_t\,\dev_t^i\,\dt - \sk{\K_t \dI_t}\,\dev_t^i + O(\epsilon^2), \label{eq:d_dev_new} 
\end{align}
where $A_t = -\Omega_t$. 

Define $\Gamma_t^i:= \dev_t^i\,\dev_t^{i\,T}$. Using the \ito's lemma, 
\begin{align*}
 \ud \Gamma_t^i = & ~\ud \dev_t^i\,\dev_t^{i\,T} + \dev_t^i\,\ud (\dev_t^{i\,T}) + \ud \dev_t^i\,\ud(\dev_t^{i\,T}) \\
 = & ~\bpar{A_t\,\dt - \sk{\K_t \dI_t}} \Gamma_t^i + \Gamma_t^i \bpar{A_t\,\dt - \sk{\K_t \dI_t}}^T 
     + \sigmaB\,\bpar{\dB_t^i\,\dev_t^{i\,T} + \dev_t^i\,\dB_t^{i\,T}} + \sigmaB^2\,I\,\dt \\
   & ~-\frac{1}{2}\,\bpar{\K_t\,H_t\,\Gamma_t^i + \Gamma_t^i\,H_t^T\,\K_t^T }\,\dt + O(\epsilon^2). 
\end{align*}
By definition, $\cov_t = \pi_t(\Gamma_t^i)$. 
Taking the conditional expectation on both sides and using the formula of $\K_t$, 
\begin{align*}
 \ud\cov_t = & ~\bpar{A_t\,\dt - \sk{\K_t \dI_t}} \cov_t + \cov_t \bpar{A_t\,\dt - \sk{\K_t \dI_t}}^T 
 + \sigmaB^2 \, I \,\dt - \frac{1}{\sigmaW^2}\,\cov_t H_t^T H_t \cov_t \, \dt,
\end{align*}
where $O(\epsilon^2)$ terms have been ignored.

\bibliographystyle{plain}
\bibliography{TAC_FPF_LG_abbr}

\end{document}